\documentclass[a4paper,
10pt
]{article}

\usepackage[english]{babel}

\usepackage[a4paper,left=3cm,right=3cm,top=2.5cm,bottom=2.5cm]{geometry}

\usepackage{hyperref}
\usepackage{bm}
\usepackage{authblk}
\usepackage{siunitx}
\usepackage{graphicx}
\usepackage{amsmath}
\graphicspath{ {./images/} }
\usepackage[table]{xcolor}
\usepackage{hhline}
\definecolor{lightgray}{gray}{0.9}
\usepackage{tikz}
\usepackage{subcaption}
\usetikzlibrary{positioning}
\usepackage[draft,inline,marginclue]{fixme}
\usepackage[numbers,sort&compress]{natbib}
\usepackage{subcaption}

\usepackage{array}
\newcolumntype{P}[1]{>{\raggedright\arraybackslash}p{#1}}

\FXRegisterAuthor{gs}{adt}{\color{red}GS}

\title{Optimized Bayesian Framework for Inverse Heat Transfer Problems Using Reduced Order Methods}

\author[1]{Kabir ~Bakhshaei\footnote{kbakhsha@sissa.it}}
\author[2]{Umberto Emil ~Morelli\footnote{umbertoe.morelli@gmail.com }}
\author[1]{Giovanni ~Stabile\footnote{giovanni.stabile@uniurb.it:  Corresponding author}}
\author[1]{Gianluigi ~Rozza\footnote{grozza@sissa.it}}

\affil[1]{Mathematics Area, mathLab, SISSA, via Bonomea 265, I-34136 Trieste, Italy}
\affil[2]{Universidade de Santiago de Compostela, Santiago de Compostela, Spain}
\affil[3]{Department of Pure and Applied Sciences, Informatics, and Mathematics Section, the University of
Urbino Carlo Bo, Piazza della Repubblica, 13, I-61029 Urbino, Italy}

\date{\today} 

\begin{document}
\maketitle
\listoffixmes
\begin{abstract}
A stochastic inverse heat transfer problem is formulated to infer the transient heat flux, treated as an unknown Neumann boundary condition. Therefore, an Ensemble-based Simultaneous Input and State Filtering as a Data Assimilation technique is utilized for simultaneous temperature distribution prediction and heat flux estimation. This approach is incorporated with Radial Basis Functions not only to lessen the size of unknown inputs but also to mitigate the computational burden of this technique. The procedure applies to the specific case of a mold used in Continuous Casting machinery, and it is based on the sequential availability of temperature provided by thermocouples inside the mold. Our research represents a significant contribution to achieving probabilistic boundary condition estimation in real-time handling with noisy measurements and errors in the model. We additionally demonstrate the procedure's dependence on some hyperparameters that are not documented in the existing literature. Accurate real-time prediction of the heat flux is imperative for the smooth operation of Continuous Casting machinery at the boundary region where the Continuous Casting mold and the molten steel meet which is not also physically measurable. Thus, this paves the way for efficient real-time monitoring and control, which is critical for preventing caster shutdowns.\\



\end{abstract}

\section{Introduction}


The evaluation of boundary heat fluxes (HF) is a critical aspect in addressing heat transfer problems, playing a significant role in the optimization\cite{rana2022heat,thanikodi2020teaching}, design \cite{AHMED2018129,khor2021influence}, and protection \cite{huang2019coupled,thornton2021analysis} of various thermal and engineering systems. Precise determination of HF \cite{indelicato2021efficient,he2020application, ren2022flight, wei2022effects} is crucial in order to guarantee the efficiency and safety of these systems, as it impacts crucial determinations regarding the selection of materials, control of processes, and environmental factors.\\

\begin{figure*}[h]
    \centering
    
    \begin{subfigure}{0.48\textwidth}
    \includegraphics[width=\linewidth]{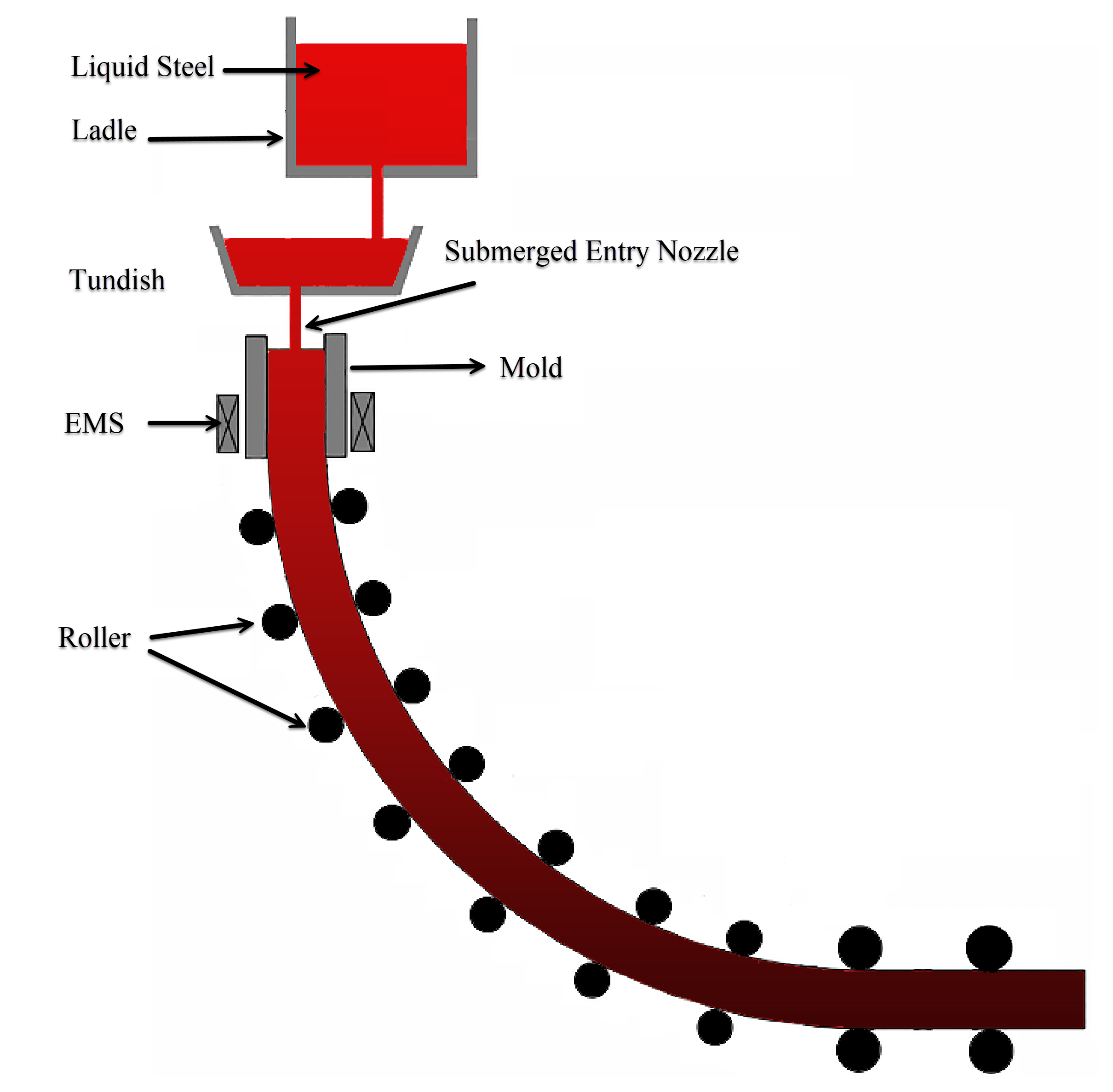} 
        \caption{}
    \end{subfigure}
    \hfill
    \begin{subfigure}{0.5\textwidth}
        \includegraphics[width=\linewidth]{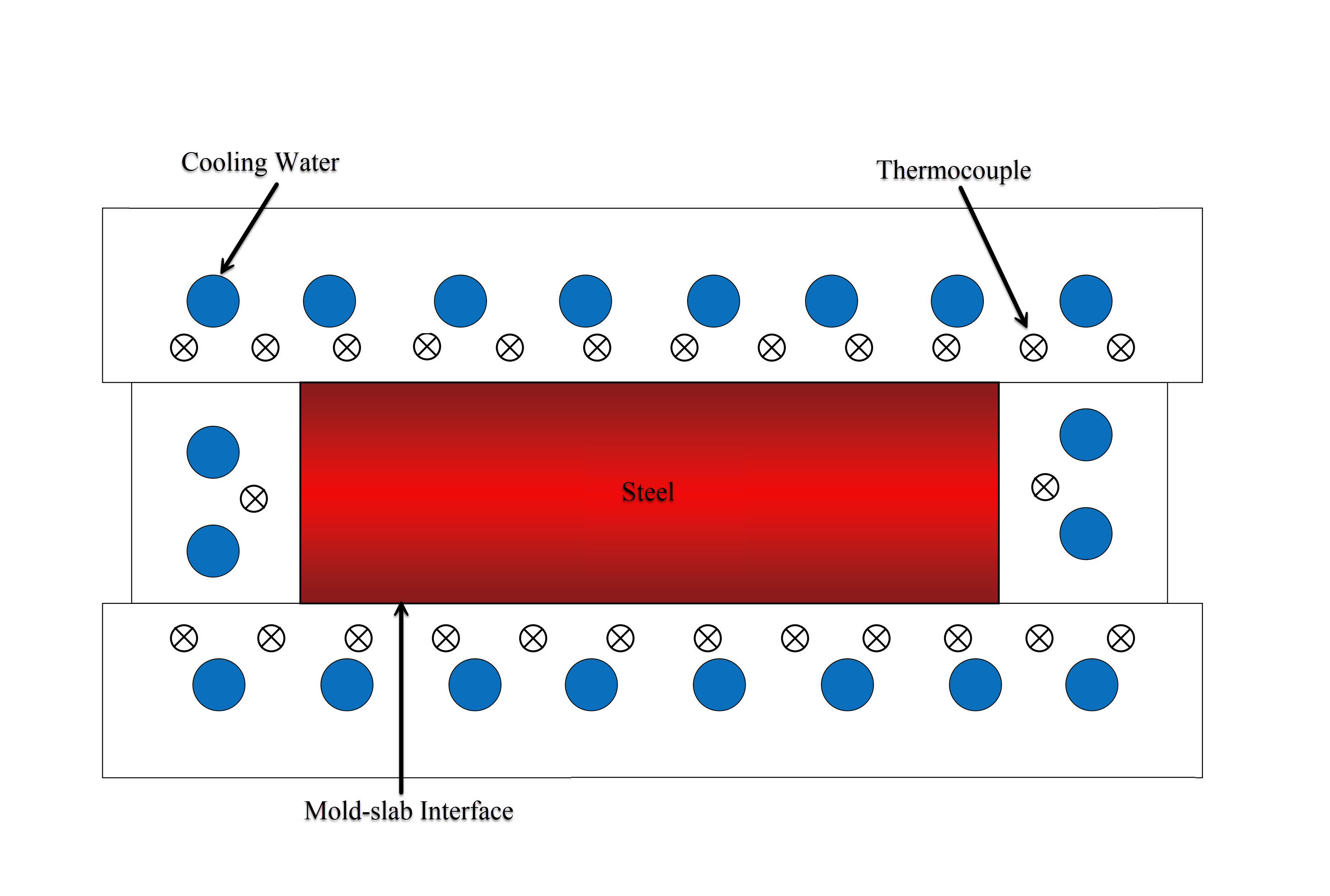}
        \caption{}
    \end{subfigure}
    \caption{ a) Schematic of CC along with b) cross-sectional view of the mold, with the casting direction oriented perpendicular to the image.}
    \label{fig:Mold-slap}
\end{figure*}

One specific problem that is the subject of our attention is predicting HF in Continuous Casting (CC) machinery, which is currently the most widely employed method for steel production. 
In the CC process, the thin slab mold, composed of copper as depicted in Figure \ref{fig:Mold-slap}, plays a crucial role in cooling the steel until it develops a thin and adequately cooled solid skin. Within the mold, a continuous circulation of running water is employed to facilitate the cooling process, while additional sensors or measurements are incorporated to oversee the temperature. Since it is impossible to put sensors directly inside the molten steel to measure its temperature, operators in CC machinery require real-time assessment of the HF at the interface between the molten steel and the CC mold in order to effectively monitor the mold's thermal behavior and avert critical issues that may necessitate the forcible shutdown of the caster \cite{Morelli2023-qw,Morelli2021-dd}. \\

Inverse Heat Transfer Problem (IHTP) in which transient HF in the interface between mold and molten steel can be considered as an unknown Neumann boundary condition is solved in a deterministic and stochastic setting. Numerous deterministic approaches have been devised \cite{BOZZOLI2014352,OKAMOTO20074409,YANG20101228,WOODBURY201331, Cabeza2005235, CHANTASIRIWAN20013823, ROUQUETTE2007128, MIERZWICZAK20102035, DOU2009728} so as to enhance the reliability of predictions in inverse heat conduction problems since this problem is frequently characterized as ill-posed due to the inherent non-uniqueness and instability of estimations based on observed temperature data. In our previous works \cite{Morelli2023-qw,Morelli2021-dd}, novel methodologies in a deterministic IHTP setting were developed to obtain steady and unsteady HF for CC machineries in the presence of noise. These innovative approaches employ an offline-online decomposition framework, which efficiently divides the computational burden between the online and offline phases, in which the offline phase is computationally expensive since it depends on mesh and time steps for solving a full-order models (FOM) but needs to be run once without relying on measurements. In addition, the sole requirement for the computationally efficient online phase to be executed exactly when the caster is working is the acquisition of thermocouple measurements. 
Thus, the proposed approaches of our previous studies take advantage of incorporating Truncated Singular Value Decomposition (TSVD)  \cite{Hansen2010-xm} regularization technique so as to overcome the ill-posedness of IHTPs\cite{Weber1981-lu} and alleviate noise propagation since small uncertainties in measurement devices can result in significant errors in determining the HF.\\

On the other hand, the Bayesian approach \cite{{WANG20043927,WANG200515, LIU20082457, PARTHASARATHY20082313}} is considered a broadly applicable and adaptable technique within stochastic methods to solve IHTPs that is inherently capable of effectively lessening the ill-posedness of IHTPs by combining the likelihood function to model the uncertainty of the data and the prior distribution to obtain posterior distributions of the unknown inputs to quantify uncertainty in the estimated parameters\cite{Cao2022-hc}. Unlike the deterministic IHTP approach in which the HF is obtained by a mean value, this quantity in the Bayesian setting is acquired as a Gaussian distribution. It is noteworthy that IHTP solutions utilizing Bayesian frameworks can incur extremely high computational costs, despite their considerable power \cite{Ramos2022-yg}. Some studies have been conducted to simultaneously determine unknown inputs along with states of linear and non-linear dynamic systems based on partial or noisy measurements. For instance, it is proved \cite{Oka2017-nr} that particle filter for predicting heat transfer coefficient and thermal conductivity in the casting process faces limitations due to sample degeneracy\cite{Fan2017-ze}, whereas, Oka et al \cite{Oka2020-yo} demonstrates that the effectiveness of EnKF in accurately estimating these two input parameters at the same time. To predict steam temperatures and heat transfer coefficients in a transient heat transfer problem involving a steam header, Khan et al \cite{Khan2021-kc} suggest a Bayesian framework using the Markov chain Monte Carlo method. Furthermore, Ramos et al \cite{Ramos2022-yg} presents a strategy for predicting the thermal conductivity and specific heat of a metal slab as changed with the temperature at the same time with the aid of Bayesian statistics and a methodology employing a transient experimental technique in three dimensions. However, none of them apply any method to improve computational efficiency. In the work of \cite{Berger2017-ah}, the thermal diffusivity coefficient of a three-dimensional, two-region heat conduction problem is predicted by employing Proper Generalized Decomposition (PGD) model reduction in conjunction with the Bayesian approach for IHTP. More examples of Bayesian method being applied to heat transfer problems are provided in \cite{kumar2018bayesian, BERGER2016327, orlande2014accelerated, 10.1115/1.4006487, GNANASEKARAN20113060, mota2010bayesian, gnanasekaran2010inexpensive, wang2004hierarchical, WANG20043927}.\\

For large-size non-linear dynamics systems, Fang et al \cite{Fang2017-tx} proposes an Ensemble-based Simultaneous Input and State Filtering (EnSISF) method with direct feed-through in a Bayesian statistical framework for estimation of wind speed and fire perimeter as unknown input and state variables, respectively, from output measurements. There are also many scenarios where employing EnSISF becomes necessary. When monitoring industrial systems that are impacted by unknown disturbances, for instance, both the operational status as a state and the disturbance itself as an input must be predicted\cite{she2004disturbance}. When tracking a maneuvering target, it is customary for the tracker to predict the target's state, including velocity and position, in addition to the input, which may consist of acceleration \cite{10.1117/12.478535}. Since Ensemble Kalman Filter (EnKF) \cite{Wackernagel2010-hy} as one of the most popular Data Assimilation (DA) techniques is primarily designed to determine the state of dynamic systems so as to deal with uncertainty and nonlinearity, in the current study, an extension of the Monte Carlo version of the Kalman Filter called EnSISF \cite{Fang2017-tx} is exploited to merge the observational data with a model to predict simultaneously not only temperature distribution as states but also an unsteady mold-slab HF as an unknown input based on the sequential availability of temperature provided by thermocouples located slightly inside the mold domain. It handles problems with a significant quantity of variables and is especially effective for highly nonlinear and high-order models. Radial Basis Functions (RBFs) are also utilized and incorporated with EnSISF to lessen the size of unknown parameters and improve the computational efficiency of the proposed approach.\\


Our research makes a significant contribution to the existing literature by applying an ensemble-based Bayesian inversion incorporated with RBFs to simultaneously determine the state variable and the amount of heat subtracted from the molten steel, which can be seen as an unknown transient boundary mold-slab HF. We utilize RBFs to reduce the size of parameters and improve computational efficiency. This method enables us to derive a probability distribution of the input when process and measurement noise are present. We also reveal the dependency of the procedure with respect to some hyperparameters like the number of seeds, shape parameter, time step, observation span, and shifting mean of the prior weights of the RBFs and their prior covariance which can not be found in the literature. It also proves that EnSISF incorporating Multiquadric is more accurate and requires less computational cost compared with integrating the approach with the Gaussian kernel.
%
%

%

The present paper is structured in the following manner. In Section \ref{sec:EnSISSF}: we explore the EnSISF incorporating RBFs concerning nonlinear systems with direct feed-through. Following this, Section \ref{sec:Dynamic Model} illustrates a CC mold that is utilized to apply the proposed algorithm. Section \ref{sec:res} studies the influence of certain hyperparameters and how they affect the accuracy of the method to estimate unknown transient HF of the mold followed by Section \ref{sec:conc} that provides a summary of our finding and offers some future research.

\section{Ensemble-based Simultaneous Input and State Filtering Algorithm with Radial Basis Functions Integration}\label{sec:EnSISSF}


In this study, the EnSISSF method with direct feedthrough (EnSISF-wDF)\cite{Fang2017-tx}, as a DA technique is incorporated with RBFs and is applied to predict input parameters and systems states with respect to the output measurements. The general form of a nonlinear dynamic system with direct input-to-output feed-through is shown below.

\begin{center}
\begin{flalign}
& \left\{
\begin{aligned}
\bm{\chi}_{k+1}=\bm{\Lambda}\left(\bm{\chi}_{k}, \bm{u}_{k}\right) +\bm{w
}_{k}, \\
\bm{y}_{k}=\bm{\Theta}\left(\bm{\chi}_{k}, \bm{u}_{k}\right) +\bm{v
}_{k}.
\end{aligned}
\right.
\label{eq:SISSFwdf}
\end{flalign}
\end{center}

At time step $k+1$, the unknown system state vector $\bm{\chi}_{k+1}$ is a nonlinear state transition function $\bm{\Lambda}$ of the previous state $\bm{\chi}_{k}$, and unknown system input vector $\bm{u}_{k}$. Moreover, $\bm{y}_{k}$ expresses a known system measurement vector, while $\bm{\Theta}$ is a nonlinear measurement function. $\bm{w}$ and $\bm{v}$ are uncorrelated zero-mean white Gaussian process and measurement noise, respectively with the respective covariance matrices $\bm{Q}_{k}$ and $\bm{R}_{k}$. Process noise $\bm{w}$, which represents the discrepancy between the model's mathematical representation and the actual occurrences, refers to unanticipated or unknown changes in a system that are not accounted for by the model.\\
 
\begin{center}
\begin{equation}
\bm{w}_{k} = \mathcal{N}(\bm{0}, \bm{Q}_{k}).
\label{eq:ProcessNoise}
\end{equation}
\end{center}

\begin{center}
\begin{equation}
\bm{v}_{k} = \mathcal{N}(\bm{0}, \bm{R
}_{k}).
\label{eq:MeasurementsNoise}
\end{equation}
\end{center}

\subsection{RBF}\label{subsec:RBF}
Instead of $\bm{u}_{k}$, a parameterization technique named RBFs via \ref{eq:gRecon} is employed to decompose the space and time dependence so as to decrease the size of unknown inputs and alleviate the computational load associated with this method.  

\begin{center}
 \begin{equation}
 \bm{u}_{k} \sim \sum_{j=1}^{N} \hat{\bm{A}}_j(t) \Phi_j(\bm{X}).
\label{eq:gRecon}
\end{equation}
\end{center}

$N$ represents the number of RBFs, and  $\hat{\bm{A}}_j(t)$ is RBF weights. Different kernels $\Phi_j(\bm{X})$ can be utilized depicted in Equation \eqref{eq:Kernel} in which $r_j$ is a radius parameter which is an Euclidean distance from a point to the center of the kernel in the case of the Gaussian and Multiquadric Kernels. Each kernel also has a shape parameter $\eta_j$ which tells how steep is each basis function. This geometrical parameter is a hyperparameter that should be operated for each kernel to check if it is good enough.\\


\begin{equation}
\Phi_j(\bm{X}) =
\begin{cases}
    e^{-(\eta_j r_j)^2} & \text{ Gaussian Kernel}, \\
    \sqrt{1 + (\eta_j r_j)^2} & \text{ Multiquadric Kernel}.
\end{cases}
\label{eq:Kernel}
\end{equation}

For each time step k = 1, 2, ..., the proposed approach works by breaking down prediction and update into ensembles, and then combining them together. The following sections express a prediction-update technique developed for the chosen strategy along with the initial step. For the assumptions and mathematical background of the method, please refer to \cite{Fang2017-tx}.\\

\subsection{Initial}\label{subsec:Initial}

The procedure starts by assuming Gaussian prior for the temperature represented in Equation \ref{eq:GausPriorstate}, and Gaussian prior distribution of unknown weights of RBFs \ref{eq:GausPriorParameter} that follows with a mean $\bm{\hat{A}}_{j,k=0}$ and its covariance represented by $\kappa\delta_{ij}\hat{A}_{j,k=0}$ in which the Kronecker delta $\delta_{ij}$ is a function that equals 1 when $i=j$ and 0 when $i\sim j$. Moreover, the variable $\kappa$ represents the scaling factor of $\bm{\hat{A}}_{j,k=0}$ to effectively scale the elements along the diagonal of the covariance matrix as a hyperparameter.

\begin{center}
\begin{equation}
\bm{\chi}_{k=0} = \mathcal{N}(\bm{T_0}, \bm{\sigma}_{T}).
\label{eq:GausPriorstate}
\end{equation}
\end{center}

\begin{center}
\begin{equation}
\bm{A}_{k=0} = \mathcal{N}(\bm{\hat{A}}_{j,k=0}, \kappa\delta_{ij} \bm{\hat{A}}_{j,k=0}).
\label{eq:GausPriorParameter}
\end{equation}
\end{center}

\subsection{Predict}\label{subsec:Predict}

At each time step, the forecast step is based on the information from the previous time step and it is used to make a prediction at the current time instant. Then, the process model $\bm{\Lambda}$ is directly solved with unknown parameters for each member of the ensemble and adding the solution to the model error \ref{eq:ProcessNoise} which eventually results in an ensemble of solution as depicted in Equation \ref{eq:forwardProblem}. Since the boundary condition is stochastic, states are also stochastic. 

\begin{center}
\begin{equation}
\widehat{\bm{\chi}}_{k \mid k-1}^i=\bm{\Lambda}\left(\hat{\bm{\psi}}_{k-1 \mid k-1}^i\right)+\bm{w}_{k-1}^i .
\label{eq:forwardProblem}
\end{equation}
\end{center}

The joint ensemble $\widehat{\bm{\psi}}_{k \mid k-1}$ is depicted in Equation \ref{eq:builtJointEns}. The sample mean of \ref{eq:forwardProblem} is calculated via \ref{eq:forwardProblemMean} in which $S_n$ is the number of ensemble.

\begin{center}
\begin{equation}
\widehat{\bm{\chi}}_{k \mid k-1} = \frac{1}{S_n} \sum_{i=1}^{S_n} \widehat{\bm{\chi}}_{k \mid k-1}^i.
\label{eq:forwardProblemMean}
\end{equation}
\end{center}

If there are no measurements available at this time step, a joint ensemble is built as shown in Equation \ref{eq:builtJointEns}, and its mean can be calculated via Equation \ref{eq:meanjointEnsembleUpdate}. 

\begin{center}
\begin{equation}
\hat{\bm{\psi}}_{k \mid k-1}^{i}=\left[\begin{array}{c}
\hat{\bm{A}}_k^{i} \\
\hat{\bm{\chi}}_{k \mid k-1}^{i}
\end{array}\right].
\label{eq:builtJointEns}
\end{equation}
\end{center}

\subsection{Update}\label{subsec:Update}

In this step, available measurements are employed to improve the accuracy of the estimation nested by an iterative method in which the iteration number is denoted by $\beta$. Therefore, forecasted values for each member of the ensemble in the $\beta$th iteration is obtained by 

\begin{center}
\begin{equation}
\hat{\bm{y}}_{k \mid k-1}^{i, \beta}=\bm{\Theta}\left(\hat{\bm{\psi}}_{k \mid k-1}^{i, \beta}\right)+\bm{v}_k^{i, \beta}.
\label{eq:measurementFunction}
\end{equation}
\end{center}

During the update phase, when measurements become available, the covariance matrix of the observations \ref{eq:covarianceMeasurementsSigmaPoint} and cross-covariance matrix \ref{eq:crossCovarianceMeasurementsSigmaPointAndSigmaPoints} between the joint variables and the observations are computed with the aid of Equations \ref{eq:meanjointEnsembleUpdate}, and \ref{eq:meanMeasurementsEnsembleUpdate} representing the mean of the joint and observation ensembles, respectively.\\

\begin{center}
 \begin{equation}
\hat{\boldsymbol{\psi}}_{k \mid k-1}^{\beta}=\frac{1}{S_n} \sum_{i=1}^{S_n} \hat{\bm{\psi}}_{k \mid k-1}^{i, \beta}.
\label{eq:meanjointEnsembleUpdate}
\end{equation}
\end{center}

\begin{center}
\begin{equation}
\hat{\bm{y}}_{k \mid k-1}^{\beta}=\frac{1}{S_n} \sum_{i=1}^{S_n} \hat{\bm{y}}_{k \mid k-1}^{i, \beta}.
\label{eq:meanMeasurementsEnsembleUpdate}
\end{equation}
\end{center}

\begin{center}
\begin{equation}
\bm{P}_{k \mid k-1}^{\bm{y}, \beta}=\frac{1}{S_n} \sum_{i=1}^{S_n} \hat{\bm{y}}_{k \mid k-1}^{i, \beta} \hat{\bm{y}}_{k \mid k-1}^{i, \beta \top}-\hat{\bm{y}}_{k \mid k-1}^{\beta} \hat{\bm{y}}_{k \mid k-1}^{\beta \top}.
\label{eq:covarianceMeasurementsSigmaPoint}
\end{equation}
\end{center}

\begin{center}
\begin{equation}
\bm{P}_{k \mid k-1}^{\bm{\psi} \bm{y}, \beta}=\frac{1}{S_n} \sum_{i=1}^{S_n} \hat{\bm{\psi}}_{k \mid k-1}^{i, \beta} \hat{\bm{y}}_{k \mid k-1}^{i, \beta \top}-\hat{\boldsymbol{\psi}}_{k \mid k-1}^{\beta} \hat{\bm{y}}_{k \mid k-1}^{\beta \top}.
\label{eq:crossCovarianceMeasurementsSigmaPointAndSigmaPoints}
\end{equation}
\end{center}

The ensemble Kalman gain shown in Equation \ref{eq:EnsembleKalmanGain} which improves the accuracy of the joint variable estimation by including observational information into the ensemble of model states and inputs, is a crucial component of this approach.

\begin{center}
\begin{equation}
\bm{K}_E = \bm{P}_{k \mid k-1}^{\bm{\psi} \bm{y}, \beta}\left(\bm{P}_{k \mid k-1}^{\bm{y},\beta}\right)^{-1}.
\label{eq:EnsembleKalmanGain}
\end{equation}
\end{center}

Equation \ref{eq:updateJointEnsemble} shows how the ensemble of joint variables $\hat{\boldsymbol{\psi}}_{k \mid k}^{i, \beta}$ is being updated by summing the prior joint ensemble $\hat{\boldsymbol{\psi}}_{k \mid k-1}^{i, \beta}$ and multiplying the ensemble Kalman gain $\bm{K}_E$ by the difference between what is predicted by states \ref{eq:measurementFunction} and what is measured by thermocouples $\bm{y}_k$.

\begin{center}
\begin{equation}
\hat{\boldsymbol{\psi}}_{k \mid k}^{i, \beta}=\hat{\boldsymbol{\psi}}_{k \mid k-1}^{i, \beta}+\bm{K}_E\left(\bm{y}_k-\hat{\bm{y}}_{k \mid k-1}^{i, \beta}\right) .
\label{eq:updateJointEnsemble}
\end{equation}
\end{center}

Thus, the mean of the updated temperatures and the RBFs weight is extracted via \ref{eq:meanOfUpdateJointEnsemble}. As a result, the reconstructed $\bm{u}_{k}$ is computed via Equation \ref{eq:gRecon}.


\begin{center}
\begin{equation}
\hat{\bm{\psi}}_{k \mid k}^{\beta}=\frac{1}{S_n} \sum_{i=1}^{S_n} \hat{\bm{\psi}}_{k \mid k}^{i, \beta}.
\label{eq:meanOfUpdateJointEnsemble}
\end{equation}
\end{center}

This algorithm for EnSISF-wDF incorporating RBFs is summarized in Table \ref{tab:EnSISF-wDF} for the system described in \ref{eq:SISSFwdf}. This algorithm consists of two main steps, forecasting and observation. It manages ensembles that depict the conditional distributions of input parameters and system states in consideration of the output measurements, and it derives estimation of joint ensembles from ensemble sample means and covariances. The observation step is iteratively performed to enhance the accuracy of the estimation. Notably, the algorithm's characteristic features include derivative-free computation and the absence of covariance matrix propagation, making it highly efficient and well-suited for applications involving high-dimensional linear and nonlinear systems\cite{Fang2017-tx}.\\

\begin{table}[h]
\caption{The algorithm for EnSISF-wDF incorporating RBFs.}
\label{tab:EnSISF-wDF}
\centering
\begin{tabular}{|P{0.8\linewidth}|} 

\hline
\textbf{Initialization:}\\

\hhline{~|}
\rowcolor{lightgray}
\hspace{2\tabcolsep} 1: use \ref{eq:GausPriorstate} and \ref{eq:GausPriorParameter} to create joint ensemble of samples $\left\{\hat{\bm{\psi}}_{0|0}^i, i=1,2, \ldots, S_n\right\}$ via \ref{eq:builtJointEns} for $k=0$\\

\hhline{~|}
\textbf{Iteration:} \\

\hhline{~|}
\rowcolor{lightgray}
\hspace{2\tabcolsep}
2: increment $k$ by 1: $k \leftarrow k + 1$ \\
\hhline{~|}
\hspace{2\tabcolsep}
\textbf{Forecasting:}\\

\hhline{~|}
\rowcolor{lightgray}
\hspace{4\tabcolsep}
3: use \ref{eq:ProcessNoise} to create the ensemble $\left\{\bm{w}_{k-1}^i \right\}$\\

\hhline{~|}
\hspace{4\tabcolsep}
4: create forecast ensemble model $\left\{\widehat{\bm{\chi}}_{k \mid k-1}^i\right\}$ by projecting the joint ensemble $\left\{\hat{\bm{\psi}}_{k-1 \mid k-1}^i\right\}$ through \ref{eq:forwardProblem} \\

\hhline{~|}
\rowcolor{lightgray}
\hspace{4\tabcolsep}
5:  
use \ref{eq:forwardProblemMean} to calculate the sample mean $\left\{\widehat{\bm{\chi}}_{k \mid k-1}\right\}$ from $\left\{\widehat{\bm{\chi}}_{k \mid k-1}^i\right\}$ calculated in step 4\\

\hhline{~|}
\hspace{2\tabcolsep}
\textbf{Observation:}\\

\hhline{~|}
\rowcolor{lightgray}
\hspace{4\tabcolsep}
6: set $\beta=0$ to start the iterative observation\\

\hhline{~|}
\hspace{4\tabcolsep}
7: \textbf{Continuing while} $\beta \leq \beta_{\max }$ \\

\hhline{~|}
\rowcolor{lightgray}
\hspace{6\tabcolsep}
8: create weight ensemble of samples $\left\{\hat{\bm{A}}_k^{i, \beta}\right\}$ from $\mathcal{N}\left(\hat{\bm{A}}_k^{\beta-1}, \kappa\delta_{ij} \bm{\hat{A}}_{j,k=0}\right)$, when $\beta=0$, $\hat{\bm{A}}_k^{\beta-1}=\mathbf{A}_k$\\

\hspace{6\tabcolsep}
9: create joint ensemble $\left\{\hat{\bm{\psi}}_{k \mid k}^{i, \beta}\right\}$ by merging $\left\{\hat{\mathbf{A}}_k^{i, \beta}\right\}$ from step 8 and $\left\{\widehat{\bm{\chi}}_{k \mid k-1}^i\right\}$ from step 4\\

\hhline{~|} 
\rowcolor{lightgray}
\hspace{6\tabcolsep}
10: use \ref{eq:MeasurementsNoise} to create the ensemble $\left\{\hat{\bm{v}}_k^{i, \beta}\right\}$\\

\hhline{~|} 
\hspace{6\tabcolsep}
11: use \ref{eq:measurementFunction} to create the ensemble $\left\{\hat{\bm{y}}_{k \mid k}^{i, \beta}\right\}$\\

\hhline{~|}
\rowcolor{lightgray}
\hspace{6\tabcolsep}
12: use \ref{eq:meanjointEnsembleUpdate}, \ref{eq:meanMeasurementsEnsembleUpdate}, \ref{eq:covarianceMeasurementsSigmaPoint}, \ref{eq:crossCovarianceMeasurementsSigmaPointAndSigmaPoints} to create the sample covariance $\bm{P}_{k \mid k}^{\bm{y}, \beta}$ and cross covariance $\bm{P}_{k \mid k}^{\bm{\psi} \bm{y}, \beta}$\\

\hhline{~|}
\hspace{6\tabcolsep}
13: execute the ensemble update using \ref{eq:EnsembleKalmanGain}, \ref{eq:updateJointEnsemble} to derive $\left\{\hat{\bm{\psi}}_{k \mid k}^{i, \beta}\right\}$\\

\hhline{~|}
\rowcolor{lightgray}
\hspace{6\tabcolsep}
14: use \ref{eq:meanOfUpdateJointEnsemble} to calculate the mean of the sample $\left\{\hat{\boldsymbol{\psi}}_{k \mid k}^{i, \beta}\right\}$\\

\hhline{~|}
\hspace{6\tabcolsep}
15: extract $\hat{\mathbf{A}}_k^{\beta}$ from $\hat{\bm{\psi}}_{k \mid k}^{\beta}$ \\

\hhline{~|}
\rowcolor{lightgray}
\hspace{6\tabcolsep}
16: increment $\beta$ by 1: $\beta \leftarrow \beta+1$ \\

\hhline{~|}
\hspace{4\tabcolsep}
17: \textbf{Conclude while loop} \\

\hhline{~|}
\rowcolor{lightgray}
\hspace{4\tabcolsep}
18: extract $\hat{\bm{A}}_{k \mid k}$ and $\hat{\bm{\chi}}_{k \mid k}$ from $\hat{\bm{\psi}}_{k \mid k}^{\beta_{\max }}$ \\

\hhline{~|}
\hspace{2\tabcolsep}
19: continue until there are no more incoming measurements\\

\hline
\end{tabular}
\end{table}
\newpage \clearpage
\section{Application to CC Mold}\label{sec:Dynamic Model}

\begin{figure}[h]
    \centering
\includegraphics[scale=0.16]{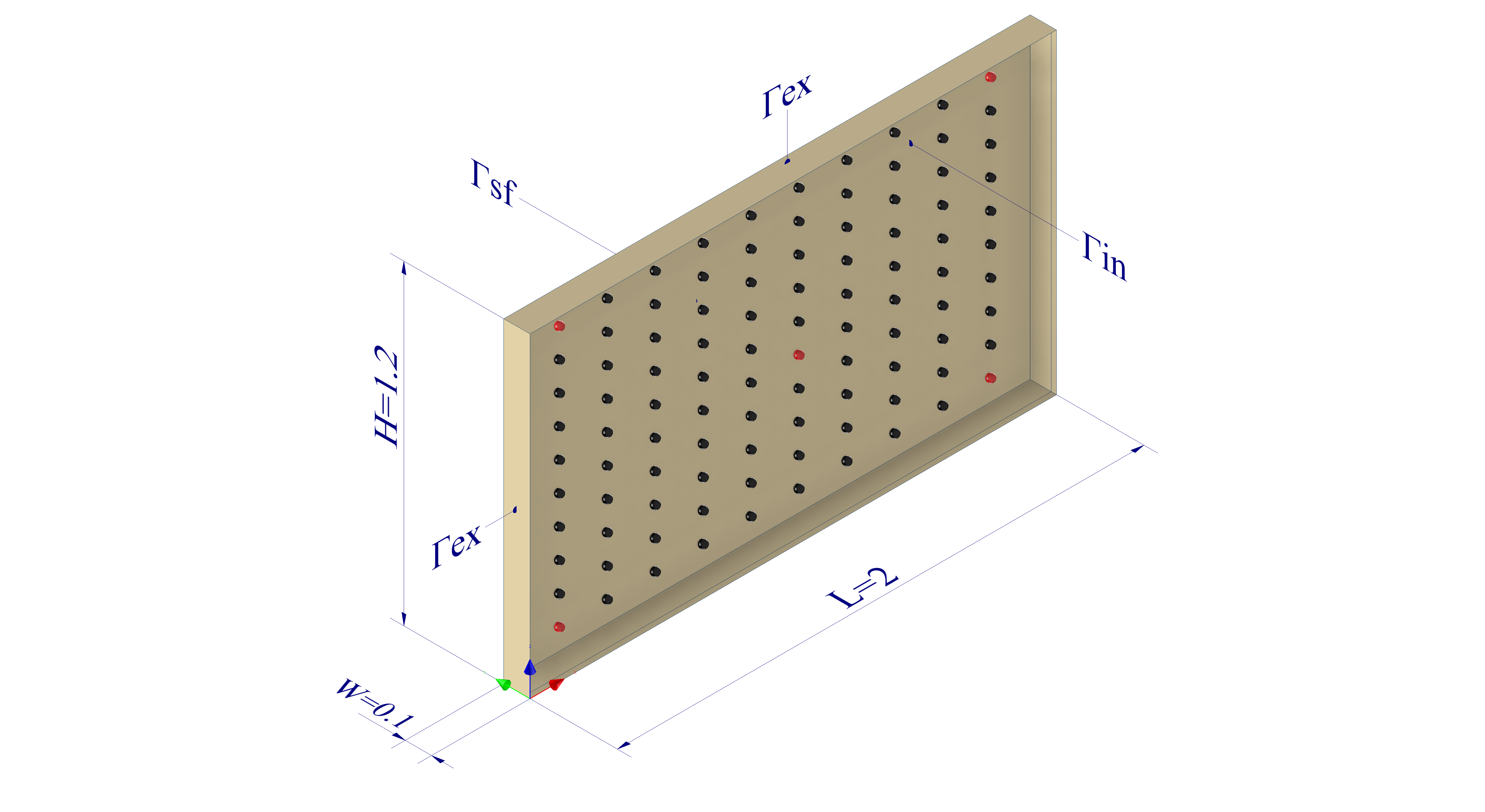}
    \caption{Transparent mold with labeled boundaries, showing centers of RBFs (red) and thermocouple locations (black).}
    \label{fig:Mold}
\end{figure}

The approach is implemented in the specific case of mold used in CC machinery as shown in Figure \ref{fig:Mold}. It is based on an unsteady 3D heat conduction model that incorporates the following assumptions described in our prior publications\cite{Morelli2023-qw,Morelli2021-dd}. The mold is treated as a solid material exhibiting homogeneity and isotropy, with minimal consideration for thermal expansion and mechanical distortion. It is assumed that the material properties remain constant throughout the analysis. The boundaries in contact with air are regarded as adiabatic, and the influence of radiation on heat transmission is intentionally omitted.\\


The resulting heat conduction model is defined by its boundary and initial conditions, which are given by Equation \eqref{eq:BoundaryCondition}.

\begin{equation}
\rho c_p \frac{\partial T}{\partial t} - k_s \Delta T = 0. \hspace{5.8cm} \text{in } \Omega \times (0, t_f]
\label{eq:3DUnsteadyHeatconduction}
\end{equation}

Where $T$ denotes the temperature distribution within the mold as a function of time $t$ and location $\bm{X}$.  The terms $\rho = 5 , \text{kg/m}^3$, $c_p = 20 , \text{J/(kg} \cdot \text{K)}$, and $k_s = 3 , \text{W/(m} \cdot \text{K)}$ represent the density, specific heat capacity, and thermal conductivity of the mold material, respectively.

\begin{flalign}
& \left\{
\begin{aligned}
-k_s\nabla T \cdot \mathbf{n} &= g(t, \mathbf{X}), && \text{on } \Gamma_{s_\text{in}} \times (0, t_f] \\
-k_s\nabla T \cdot \mathbf{n} &= 0, && \text{on } \Gamma_{s_\text{ex}} \times (0, t_f] \\
-k_s\nabla T \cdot \mathbf{n} &= h(T - T_f), && \text{on } \Gamma_{s_\text{sf}} \times (0, t_f] \\
T(\cdot, 0) &= T_0. && \text{in } \Omega \\
\end{aligned}
\right. 
\label{eq:BoundaryCondition}
\end{flalign}

Where the initial temperature for the domain is ${T}_\text{0} = 400 \text({K})$. Furthermore, $g(t, \mathbf{X})$ represents an unknown boundary condition while $h = 5.66e4 , \text{W/(m}^2 \cdot \text{K)}$ expresses the heat transfer coefficient. The Laplacian operator $\Delta$ defines the spatial variation of the temperature distribution.
$\Gamma = \Gamma_{S_\text{in}} \cup \Gamma_{S_\text{ex}} \cup \Gamma_{S_\text{f}}$
 represents the boundary of $\Omega$, and it is composed of three distinct parts: $\Gamma_{S_\text{in}}$, $\Gamma_{S_\text{ex}}$, and $\Gamma_{S_\text{sf}}$.\\

The model incorporates specific boundary conditions that vary depending on the region of the mold's surface. A convective boundary condition (also called cold side boundary) is applied to the portion of the boundary $\Gamma_{S_\text{sf}}$ that is in contact with the cooling water with a temperature of ${T}_\text{f} = 350 \text({K})$, an unknown Neumann boundary condition is imposed on the portion in contact with the molted steel known as hot side boundary $\Gamma_{S_\text{in}}$, and an adiabatic boundary condition is considered for the portion in contact with the air $\Gamma_{S_\text{ex}}$.\\

The process model $\bm{\Lambda}$ in the forward problem shown in Equation \eqref{eq:SISSFwdf} is not computationally expensive for this case study, and as a result, the FOM can be employed. The FOM simulation of the state transition function is conducted using OpenFOAM\cite{Moukalled2016-ws}, a well-known open-source software for CFD applications. OpenFOAM leverages the finite-volume (FV) method \cite{moukalled2016finite}, to discretize and solve the unsteady heat conduction model. In this approach, small control volumes (cells) are utilized to determine the computational domain in which the heat conduction equation is integrated. By using the divergence, also known as Gauss's theorem, the volume integrals are converted into surface integrals, which are further discretized as summations of fluxes at the boundary faces of each control volume. By approximating these fluxes, the heat conduction equation is converted to a discrete algebraic equation system, allowing for numerical solutions using iterative or direct methods. For Uncertainty Quantification (UQ) of the method, MUQ \cite{Parno2021} library developed at MIT is coupled with ITHACA-FV \cite{noauthor_undated-xd, Stabile2018-xk, Stabile2017-sv} which is an in-house C++ library developed by SISSA Mathlab based on OpenFOAM.\\

\newpage \clearpage
\section{Numerical results and discussion}\label{sec:res}
 
\begin{figure}[h]
    \centering
\includegraphics[scale=0.15]{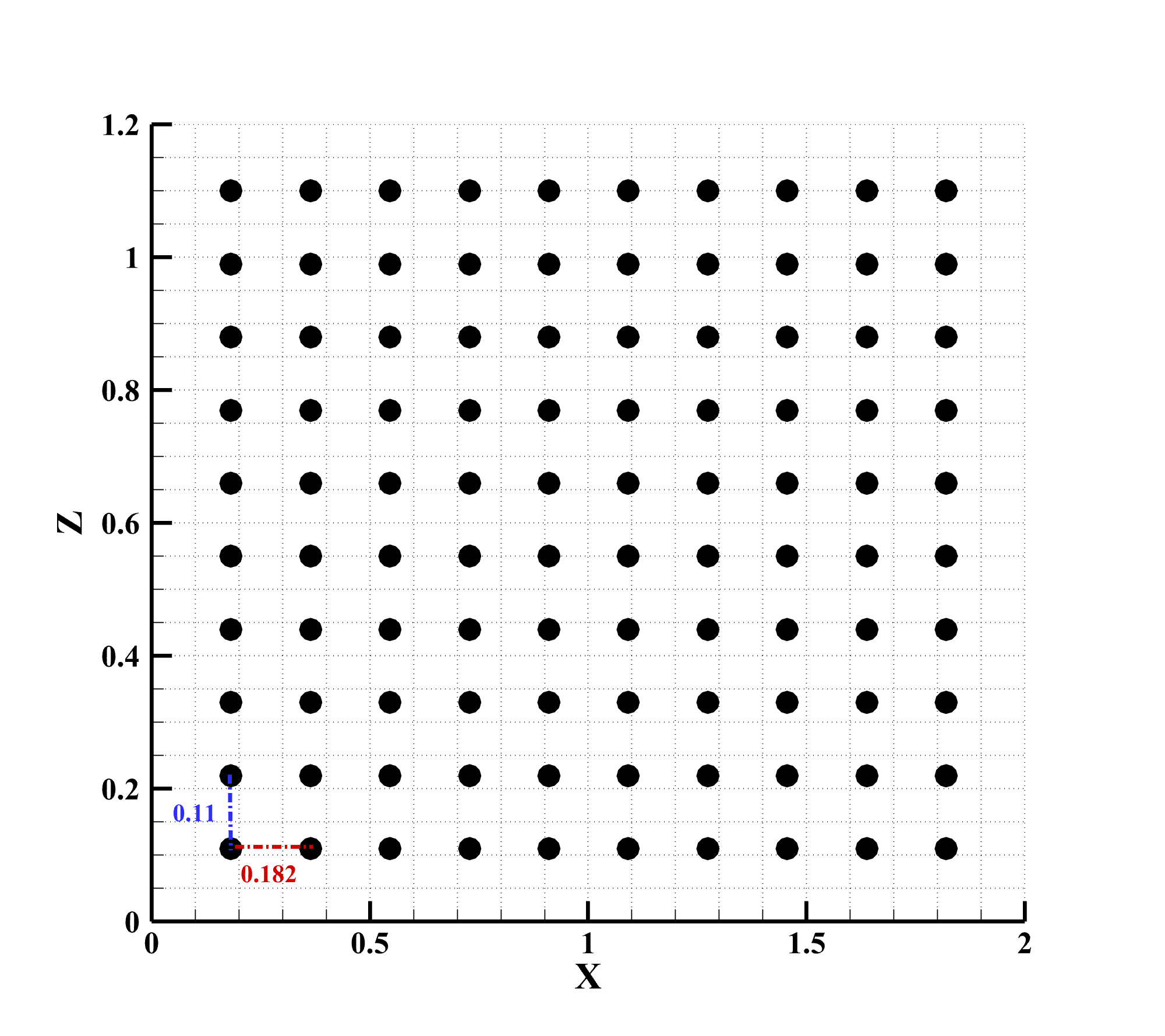}
    \caption{Position of 100 thermocouples at the plane y=0.02 m.}
    \label{fig:thermocouplespositions}
\end{figure}

This section presents the numerical results of the mold to evaluate the performance of the proposed method so as to predict and update not only states but also an unknown boundary condition represented as $g(t, \mathbf{X})$ in Equation \eqref{eq:BoundaryCondition} using sequential availability of observed temperature provided by thermocouples located at the plane $y=0.02 m$, few millimeters inward from $\Gamma = \Gamma_{S_\text{in}}$ boundary as shown in Figures \ref{fig:Mold} and \ref{fig:thermocouplespositions}. The results of a direct CFD simulation with known HF represented by $gTrue(t, \mathbf{X})$ in Equations \eqref{eq:gTrue}, \eqref{eq:ConstantforgTrue} is utilized as experimental data.

\begin{center}
\begin{flalign}
& \bm{\text{gTrue}}(\bm{X},t) = \bm{g}_1(\bm{X}) + \frac{\bm{g}_1(\bm{X})}{2} \sin\left(2\pi \cdot f_{\text{max}} \cdot \frac{t^2}{t_f}\right) + \bm{g}_2(\bm{X}) \cdot e^{-0.1 \cdot t}.
\label{eq:gTrue}
\end{flalign}
\end{center}

\begin{flalign}
& \left\{
\begin{aligned}
g_1(\bm{X}) &= -\text{k}_s \cdot (b \cdot \text{z}^2 + c), \\
g_2(\bm{X}) &= -\text{k}_s \cdot \frac{c}{1 + (\text{x} - 1)^2 + \text{z}^2} \cdot
\end{aligned}
\right.
\label{eq:ConstantforgTrue}
\end{flalign}

Where maximum Frequency and final time are $f_{\text{max}} = 0.1 \text({Hz})$ and $t_{\text{f}} = 20 \text({s})$, respectively, while $b = 200$, and $c = 300$ are constant values.\\

%
%

In both deterministic and stochastic settings, the quantity of RBFs can be arbitrary depending on the complexity of HF. However, in our previous study \cite{Morelli2023-qw,Morelli2021-dd} to solve the unknown boundary HF in a deterministic setting, the number of RBFs is the same as that of thermocouples. If the HF is complex, more basis functions should be selected. For this study, $5$ RBFs are used and their centers are defined by projecting 5 selected measurement points on the surface of $\Gamma_{S_\text{in}}$ as shown in red color in Figure \ref{fig:Mold}. Since the number of measurements is 100, the ratio of measurements to RBFs is 20. Moreover, the Euclidean distance in Equation \ref{eq:Kernel} is calculated between face centers on the $\Gamma_{S_\text{in}}$ and center of RBFs as shown in Figure \ref{fig:Mold}\\

%
%

%
    
%

%
    
%

%
\begin{figure*}[h]
    \centering
    
    \begin{subfigure}{0.6\textwidth}
    \includegraphics[width=\linewidth]{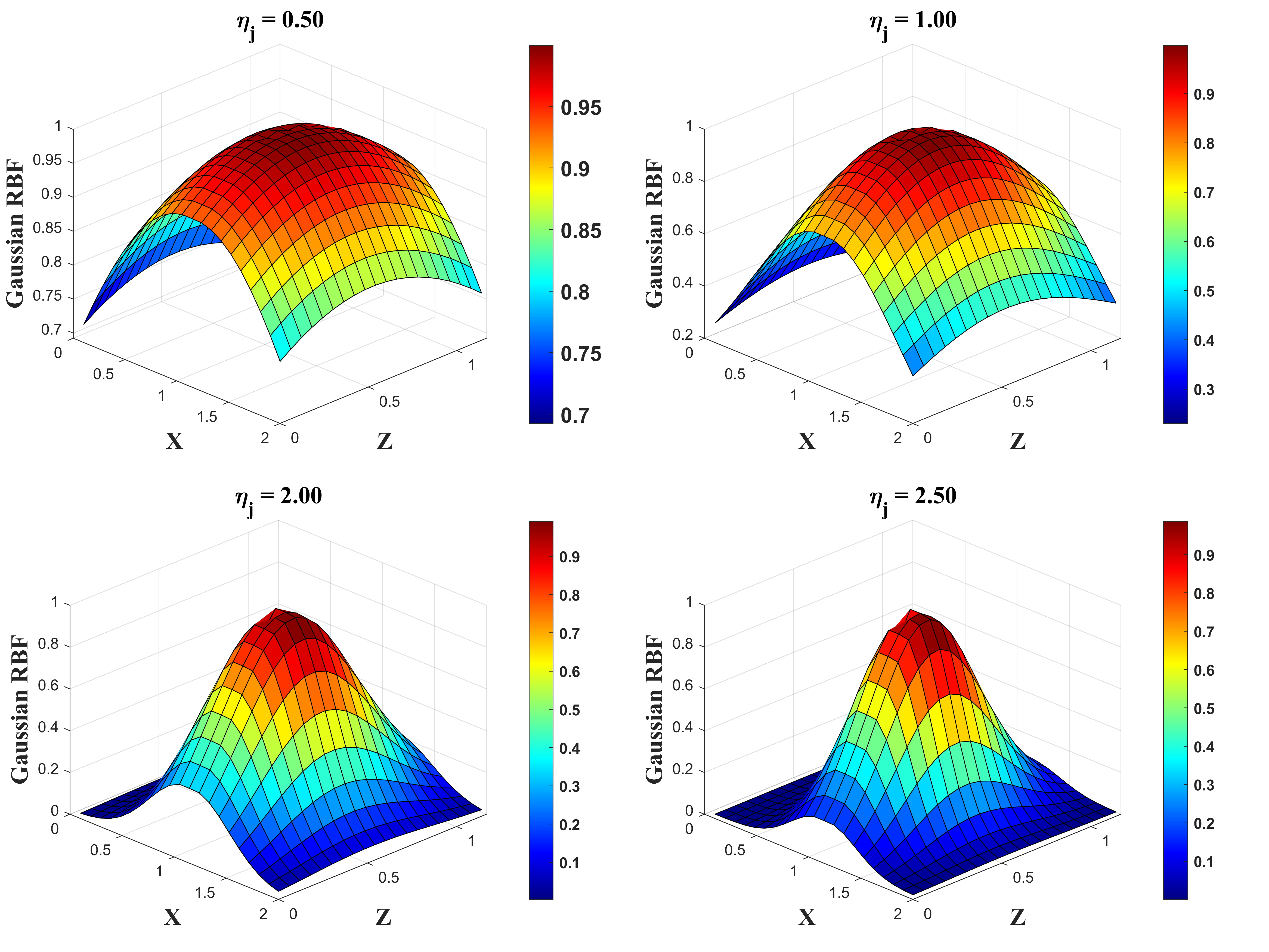} 
        \caption{}
    \end{subfigure}
    \hfill
    \begin{subfigure}{0.6\textwidth}
        \includegraphics[width=\linewidth]{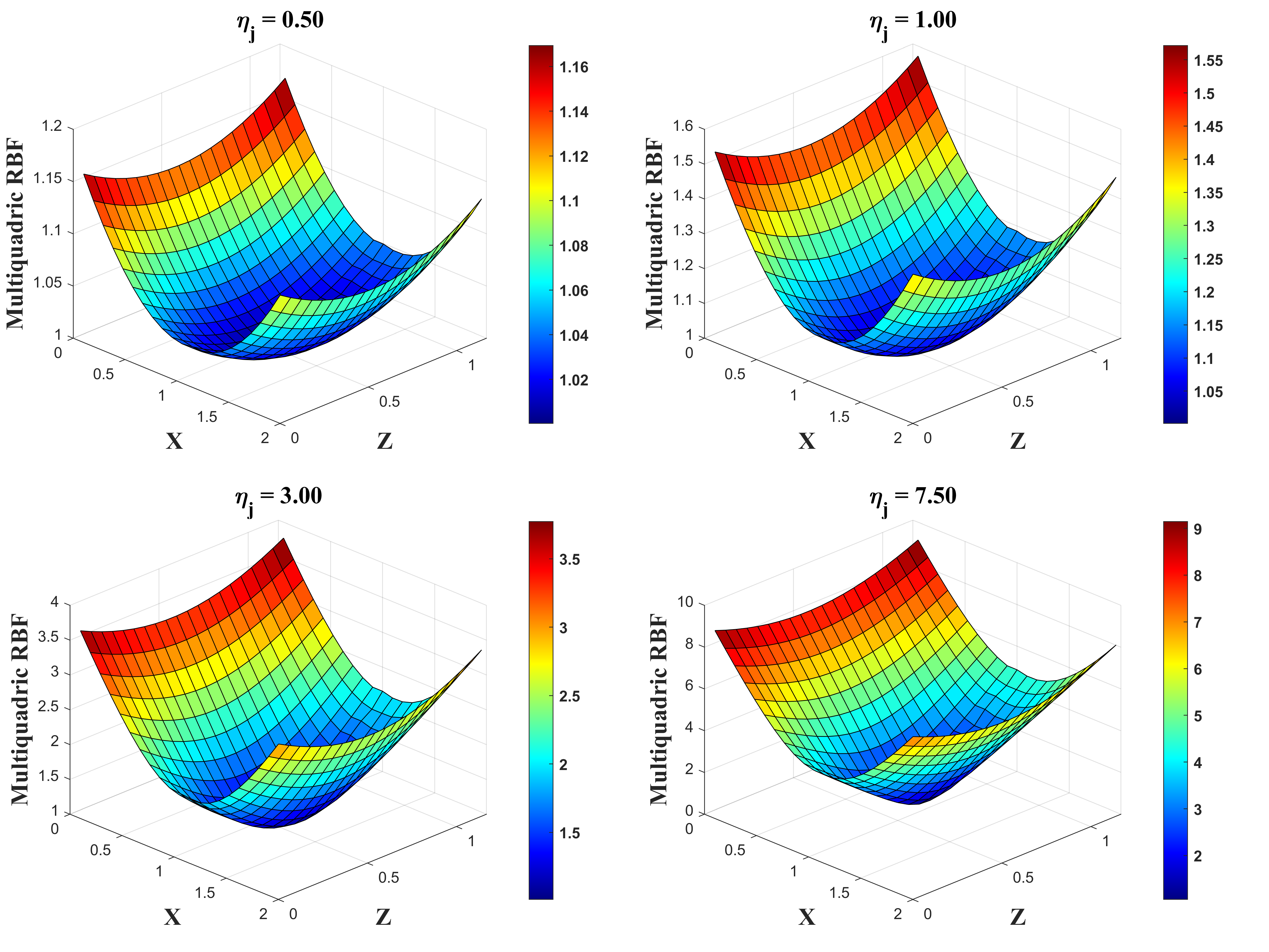}
        \caption{}
    \end{subfigure}
    \hfill
    \begin{subfigure}{0.6\textwidth}
        \includegraphics[width=\linewidth]{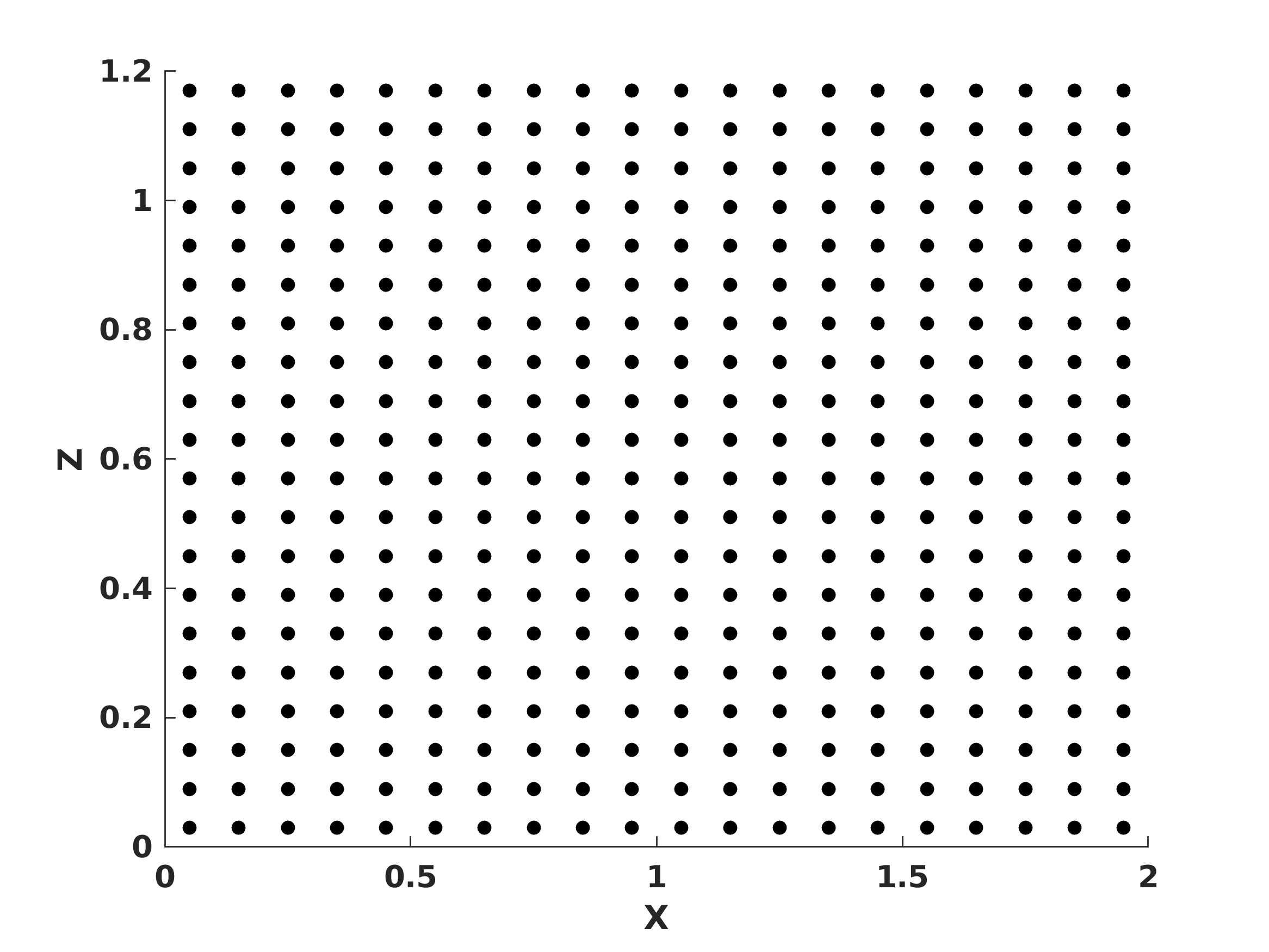}
        \caption{}
    \end{subfigure}
    \caption{RBF profiles center situated in the middle of $\Gamma_{S_\text{in}}$ for different $\eta_j$; (a) Gaussian (b) Multiquadric (c) computational nodes.}
    \label{Fig:Kernels}
\end{figure*}
%

Figure \ref{Fig:Kernels} portrays the 3D plot of Gaussian and Multiquadric RBFs in which they have the same center situated in the middle of $\Gamma_{S_\text{in}}$ shown in \ref{fig:Mold} for different shape parameter values along with computational nodes. Figure \ref{Fig:Kernels} also offers valuable insights into how variations in the shape parameter impact the shape, spread, and influence of the kernel function within the computational domain. For Gaussian, with $\eta_j = 0.5$  and Multiquadric with $\eta_j = 3$, they move towards the end of the domain, and are still different from zero indicating they are quite dispersed through the domain.\\


In order to obtain the mean of the prior weight, $\hat{\bm{A}}_{j,k=0}$, the true HF in Equation \eqref{eq:gTrue} is projected on the basis functions given by RBF with the aid of solving the following linear system of equations.\\

\begin{center}
\begin{equation}
\hat{\bm{A}}_{j,k=0} = \left(\Phi_j(\bm{X}) \cdot \Phi_j(\bm{X})^T\right)^{-1} \cdot \text{gTrue}_{k=0}(\bm{X}) \cdot
\label{fig:WeightAtT=0}
\end{equation}
\end{center}

\begin{table}[h]
\caption{Different parameters used for the mold case.}
\label{tab:parameters}
\centering
\begin{tabular}{|l l|}

\hline
\rowcolor{lightgray}
\textbf{Parameter} & \textbf{Value}\\

\hhline{~|}
$\bm{{Q}_{k}}$ & ($0.5\bm{I}$) \(K^2\)
\\
\hhline{~|}
\rowcolor{lightgray}
$\bm{{R}_{k}}$ & ($0.034\bm{I}$) \(K^2\)\\

\hhline{~|}
$\bm{\sigma}_{T}$ & ($10\bm{I}$) \(K^2\)\\

\hhline{~|}
\rowcolor{lightgray}
$\beta_{\max }$ & $1$\\

\hhline{~|}
HF probe & $(0.91, 0.00, 0.55)$ m\\

\hhline{~|}
\rowcolor{lightgray}
Temperature probe & $(0.91, 0.02, 0.55)$ m\\

\hline
\end{tabular}
\end{table}

A spatiotemporal relative error metric is utilized to distinguish the impact of hyperparameter variations on the accuracy of the proposed method. This metric compares the reconstructed HF with the actual HF across both spatial and temporal dimensions. This requires the mean relative error to be computed for both space and time. The spatiotemporal relative error values that are obtained offer significant insights sensitivity of the approach to changes in hyperparameter settings.\\

\begin{figure*}[h]
    \centering
    
    \begin{subfigure}{0.48\textwidth}
    \includegraphics[width=\linewidth]{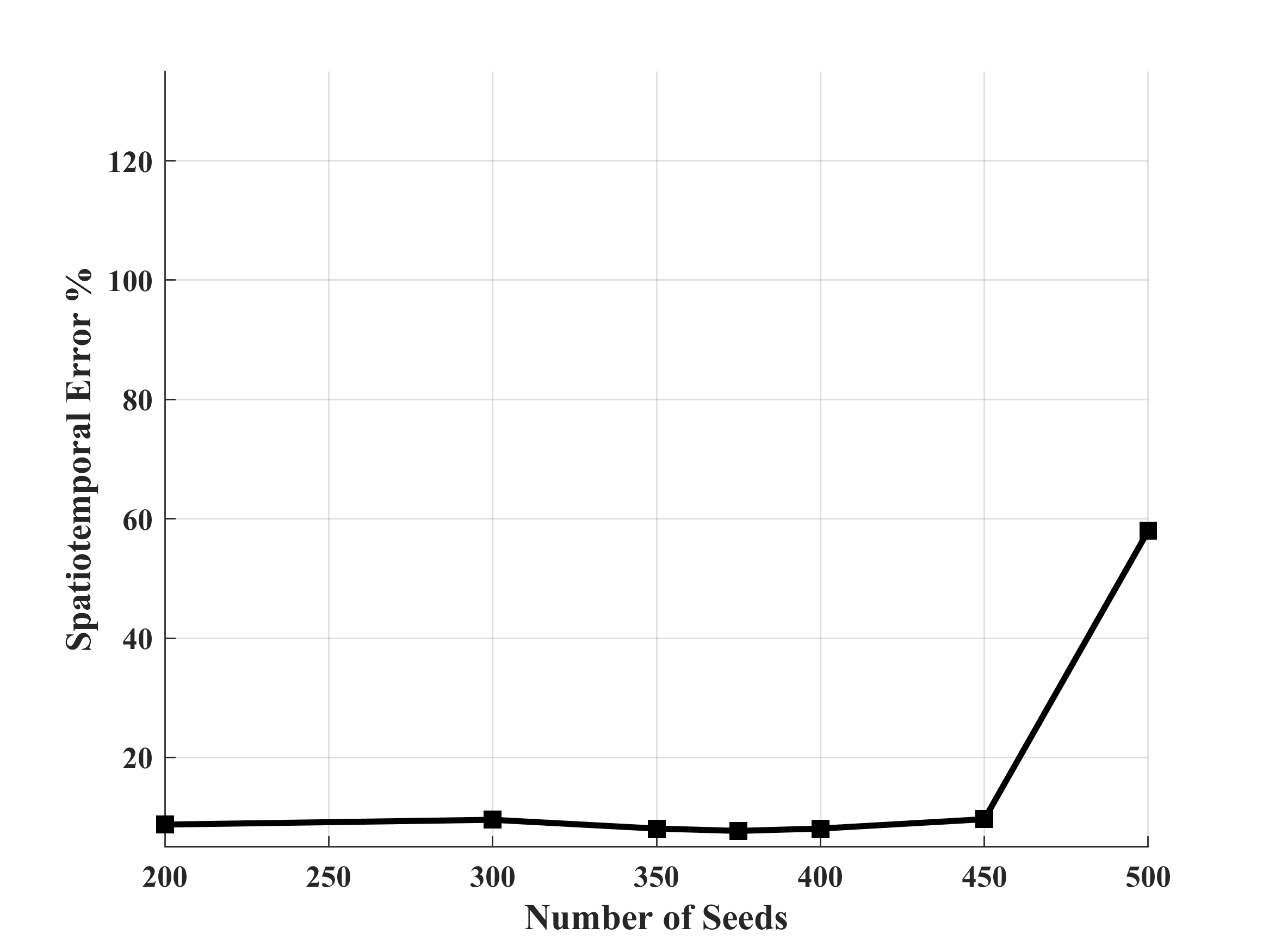} 
        \caption{}
    \end{subfigure}
    \hfill
    \begin{subfigure}{0.48\textwidth}
        \includegraphics[width=\linewidth]{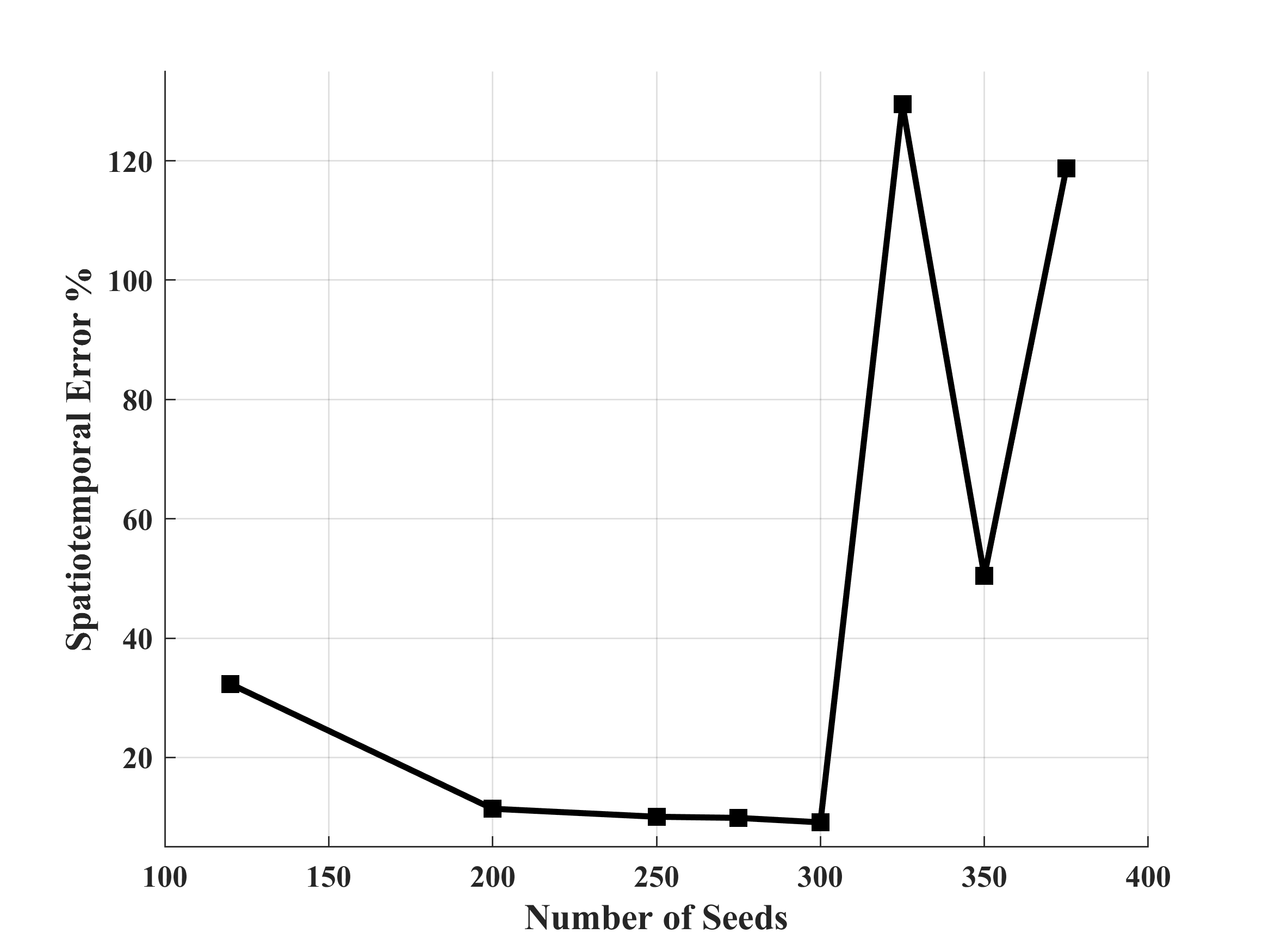}
        \caption{}
    \end{subfigure}
    \caption{Impact of the ensemble size on the spatiotemporal error; (a) Gaussian (b) Multiquadric.}
    \label{Fig:NumberOfSeeds}
\end{figure*}

For the ensemble-based method, the first hyperparameter that should be tuned is the quantity of ensembles. A nuanced relationship exists between the number of seeds and the convergence of the method. As illustrated in Figure \ref{Fig:NumberOfSeeds}, when the number of samples is increased, it leads to a reduction in the spatiotemporal error. After a specific number of samples, as this quantity increases, the errors significantly escalate as a result of the rapid rise in the condition number of the linear system that is solved in the update phase when computing the Kalman gain. When using Gaussian as a kernel, the optimal number of ensembles is observed at 375 where the spatiotemporal error is minimized at 7.7 percent, while this quantity when utilizing Multiquardic kernel decreases to 300, resulting in an error of 9.13 percent.  In the subsequent part of the paper, it is determined that the error linked to the Multiquadric kernel is reduced in comparison to the Gaussian kernel through the adjustment of additional hyperparameters as shown in Table \ref{tab:Optimalparameters}. This finding suggests that in order to reconstruct the HF for this case study the EnSSIF Multiquadric basis function incurs a lower computational cost compared with the Gaussian kernel due to its reduced sample requirement.\\

\begin{figure*}[h]
    \centering
    
    \begin{subfigure}{0.48\textwidth}
    \includegraphics[width=\linewidth]{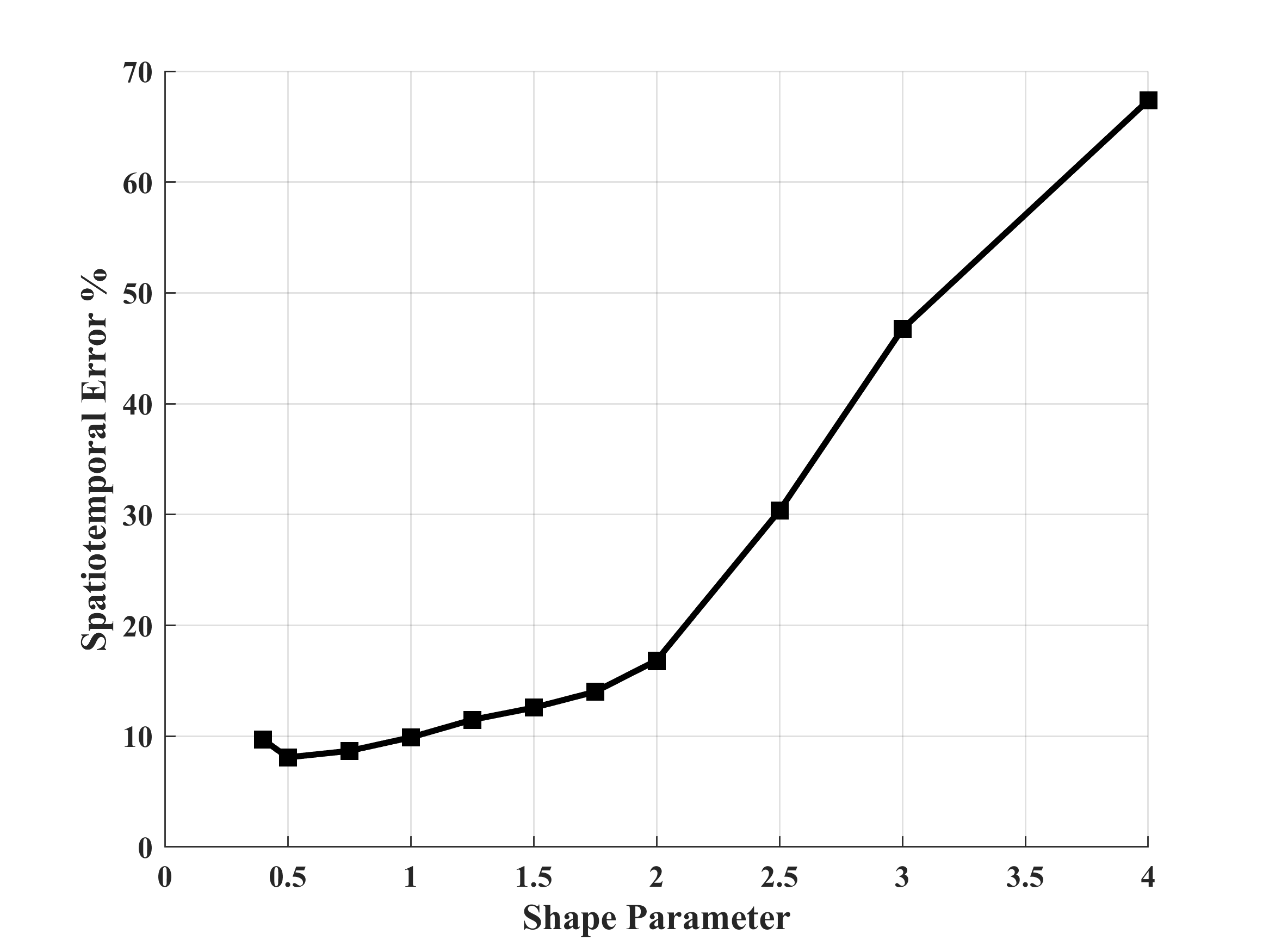} 
        \caption{}
    \end{subfigure}
    \hfill
    \begin{subfigure}{0.48\textwidth}
        \includegraphics[width=\linewidth]{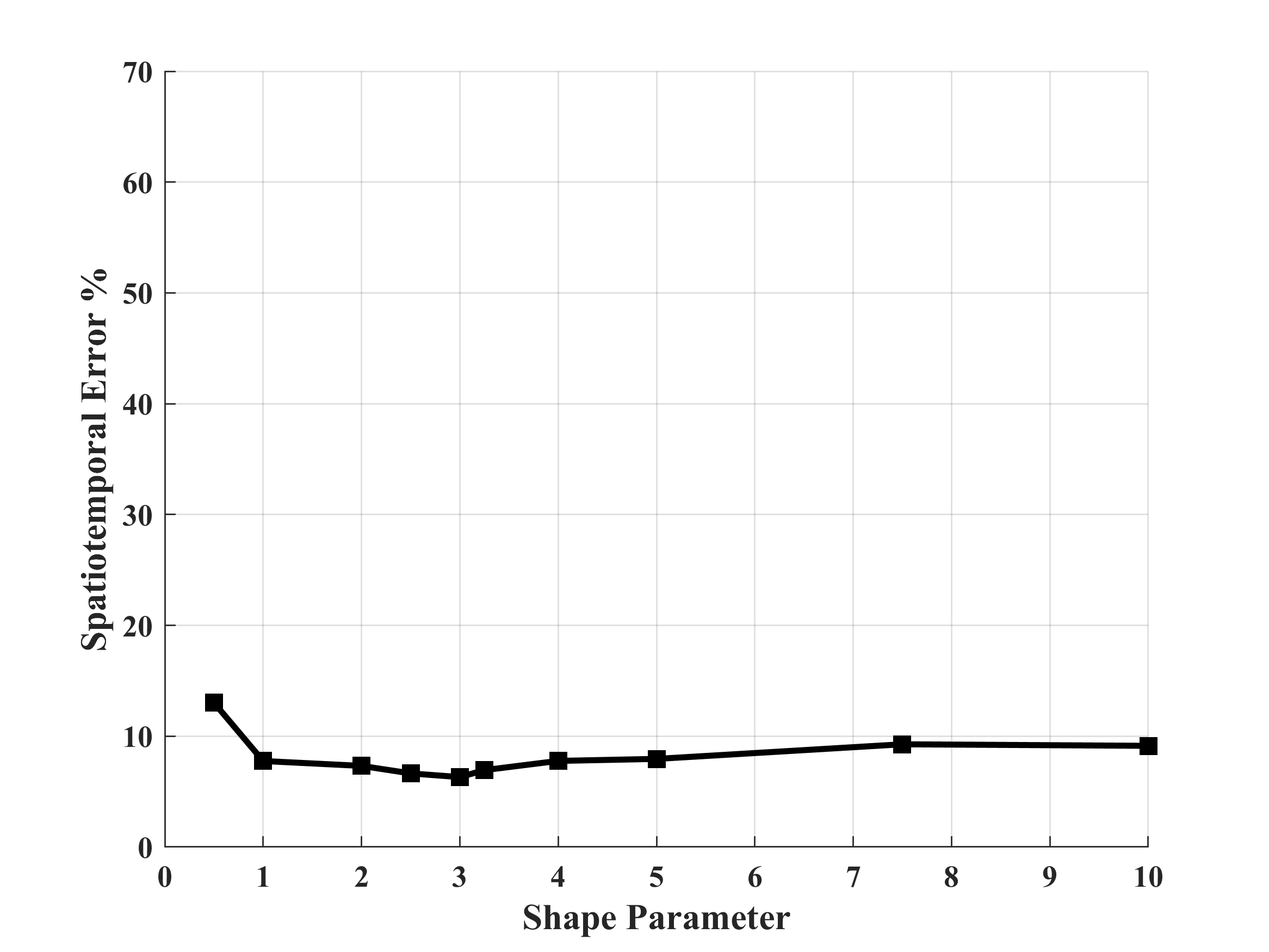}
        \caption{}
    \end{subfigure}
    \caption{Impact of shape parameter $\eta$ on the spatiotemporal error; (a) Gaussian (b) Multiquadric.}
    \label{Fig:ShapeParameter}
\end{figure*}

The geometrical parameter $\eta$ that governs the width of each kernel's distribution and the accuracy of each kernel should be appropriately adjusted. An inappropriate $\eta$ leads to an excessive amount of overlap between RBFs which eventually can undermine the accuracy and stability of the RBF interpolation. As illustrated in Figure \ref{Fig:ShapeParameter}, the optimal $\eta$ for Gaussian kernel is $0.5$, while that of the Multiquadric is $3$.\\

\begin{figure*}[h]
    \centering
    
    \begin{subfigure}{0.48\textwidth}
    \includegraphics[width=\linewidth]{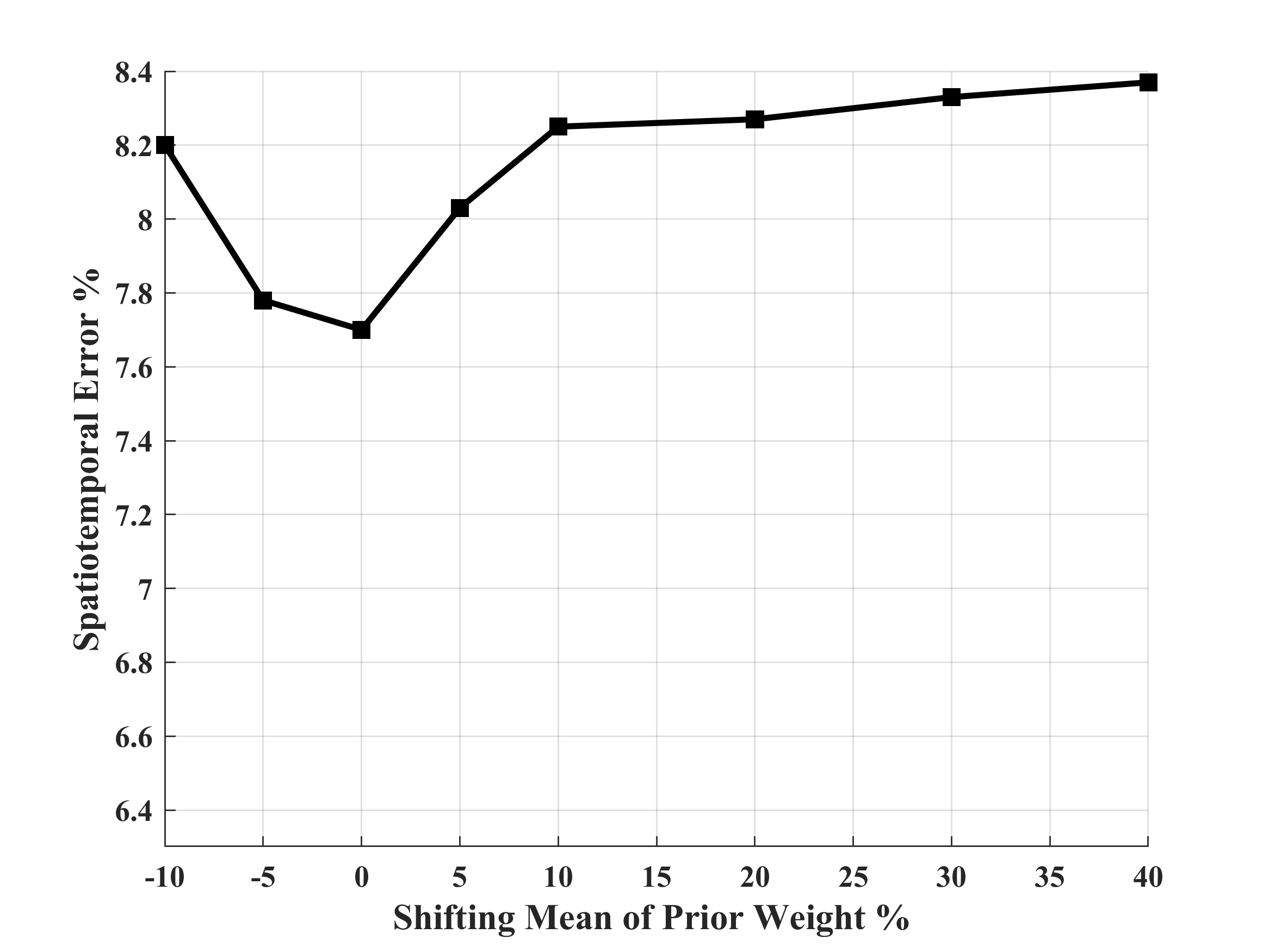} 
        \caption{}
    \end{subfigure}
    \hfill
    \begin{subfigure}{0.48\textwidth}
        \includegraphics[width=\linewidth]{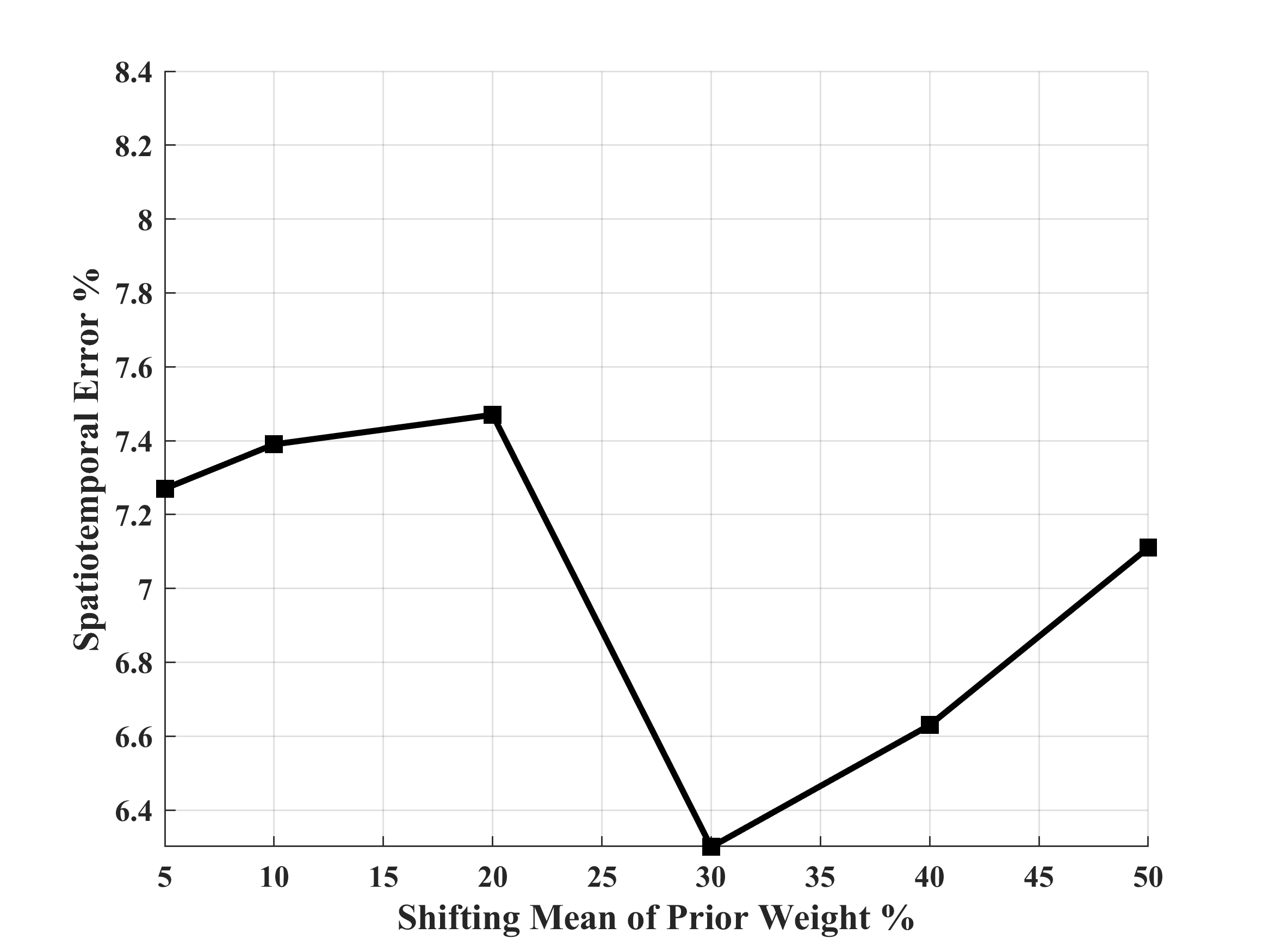}
        \caption{}
    \end{subfigure}
    \caption{Impact of the prior weight shifting on the spatiotemporal error; (a) Gaussian (b) Multiquadric.}
    \label{Fig:ShiftingMeanOfPriorWeight}
\end{figure*}

Shifting the mean of the prior weight given by the projection of the true HF onto RBFs so as not to start with a kind of exact mean is another parameter that affects the spatiotemporal error. For the Multiquadric kernel, as shown in the line chart \ref{Fig:ShiftingMeanOfPriorWeight}, the error at the initial steps consistently increases as shifting is increased, then adding some sort of shifting regularize a bit problem resulting in a reduction in error, while excessive shifting causes the error to recommence to increase again. The optimal error for Multiquadric is at 6.3 percent when 30 percent shifting is applied. Nevertheless, in the case of the Gaussian kernel, the error is optimal at 7.7 percent when there is no shifting, while from zero shifting on-wards, it experiences a rise.\\

\begin{figure*}[h]
    \centering
    
    \begin{subfigure}{0.48\textwidth}
    \includegraphics[width=\linewidth]{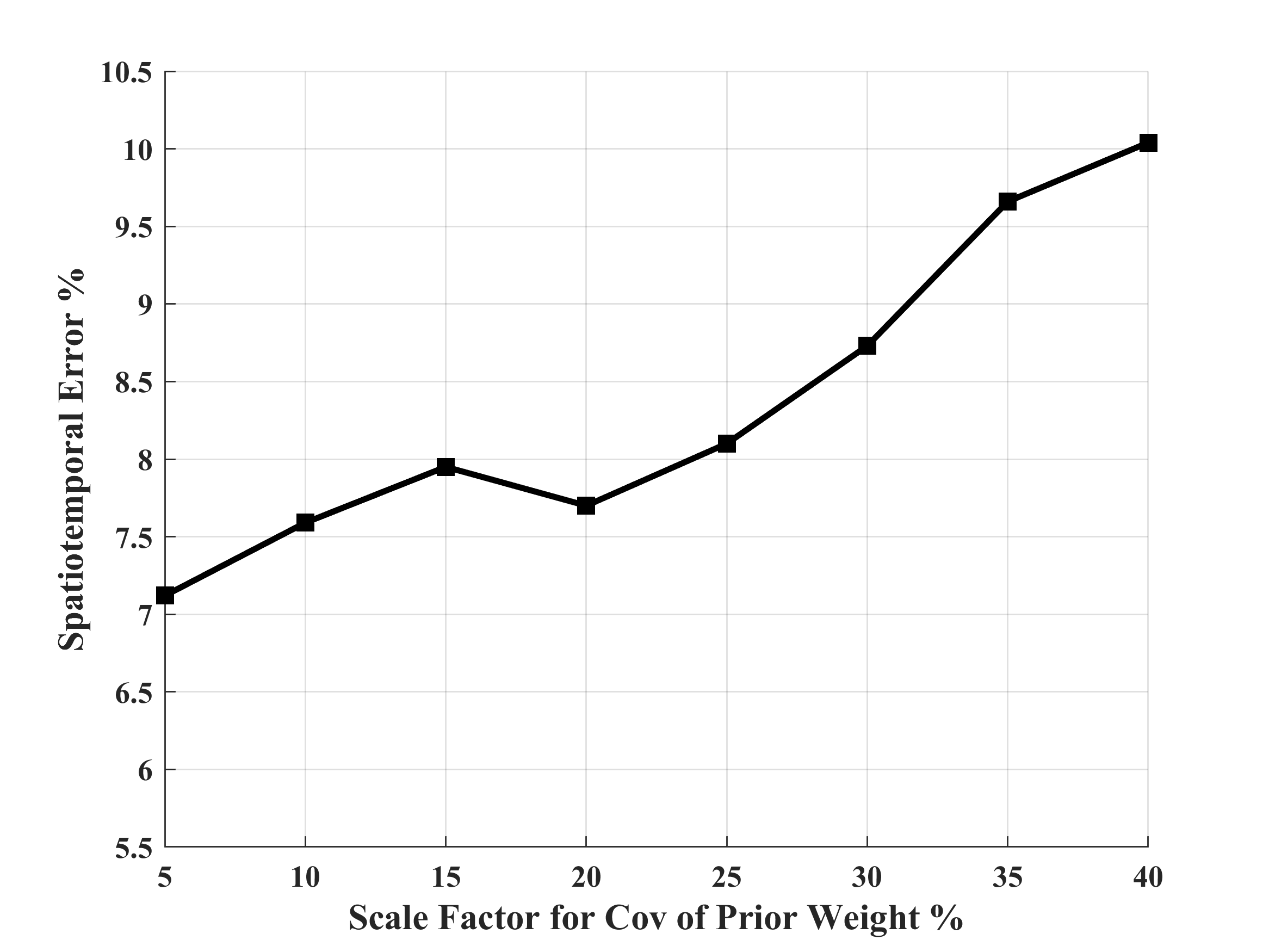} 
        \caption{}
    \end{subfigure}
    \hfill
    \begin{subfigure}{0.48\textwidth}
        \includegraphics[width=\linewidth]{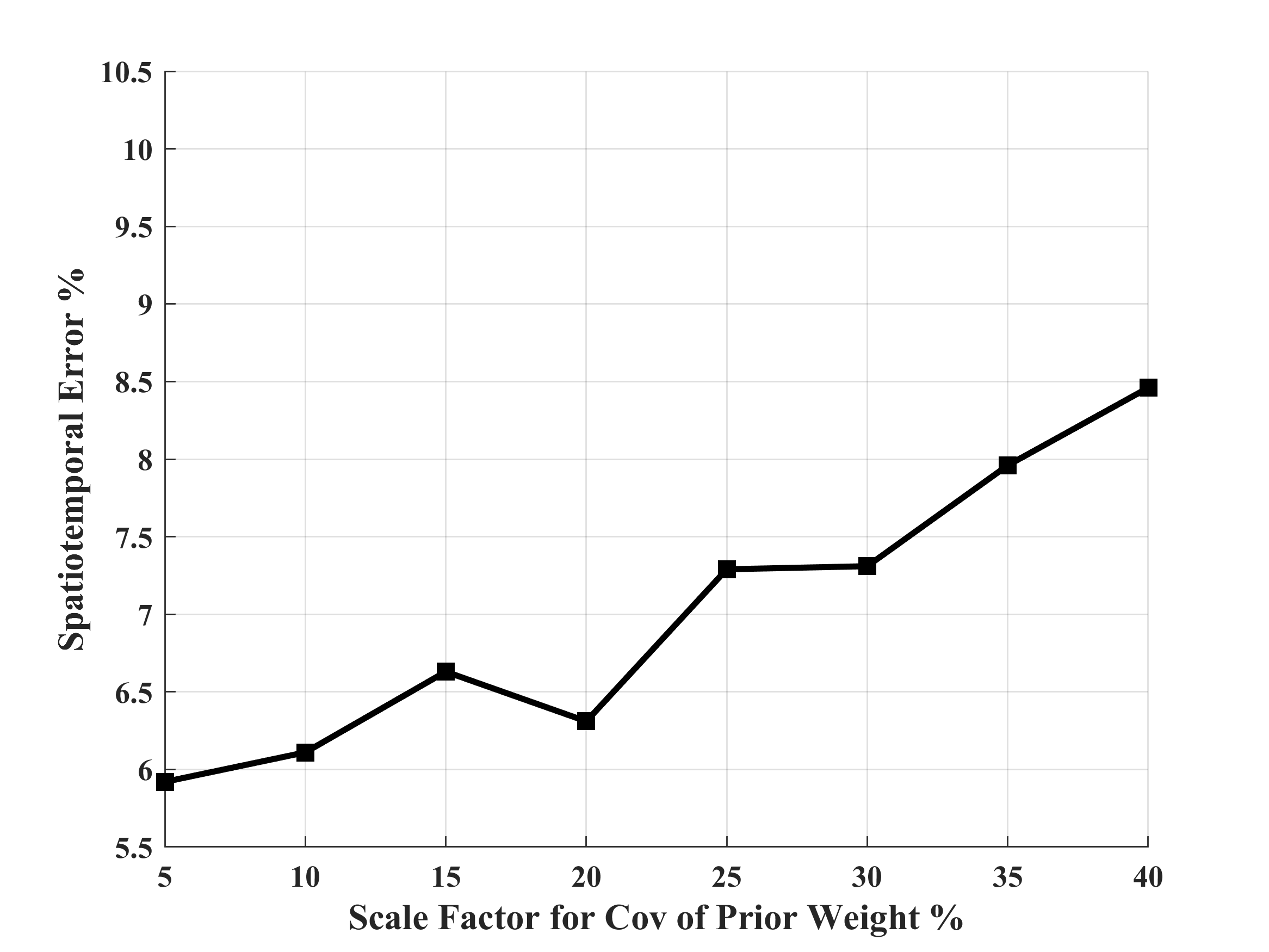}
        \caption{}
    \end{subfigure}
    \caption{Impact of $\kappa$ on the spatiotemporal error; (a) Gaussian (b) Multiquadric.}
    \label{Fig:ScaleFactorForCovOfPriorWeight}
\end{figure*}

For the covariance matrix of the prior weight as shown in Equation \ref{eq:GausPriorParameter}, $\kappa$ percent of the mean of the prior weigh is taken to see the impact of this varying scale value on the spatiotemporal error. 
The graph \ref{Fig:ScaleFactorForCovOfPriorWeight} reveals that initially, the error for the Multiquadric kernel climbs moderately reaching a peak at 6.63 percent.  Following the peak, the trend shifts, and there is a considerable fall in the error. This drop faces an optimal error at 6.31 when  $\kappa$ = 20 percent. However, the trend reverses once more and the error experiences a gradual growth.
Moreover, the error when Gaussian is employed has the same trend to that of Multiquadric kernel with the optimal error at 7.7 percent. It is observed that as Multiquadric kernel is utilized the error is less than that of Gaussian for the same $\kappa$.\\

\begin{figure*}[h]
    \centering
    
    \begin{subfigure}{0.48\textwidth}
    \includegraphics[width=\linewidth]{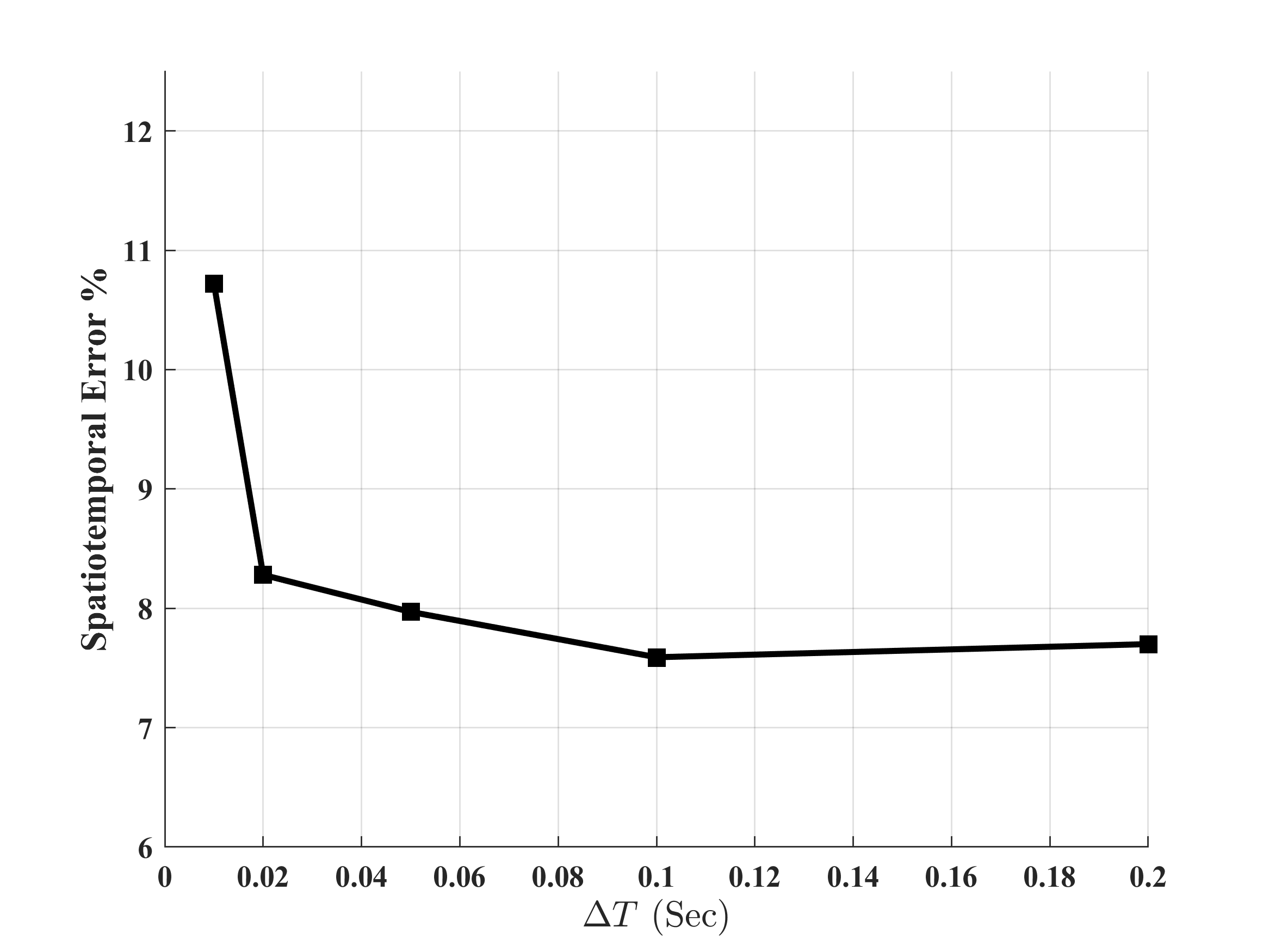} 
        \caption{}
    \end{subfigure}
    \hfill
    \begin{subfigure}{0.48\textwidth}
        \includegraphics[width=\linewidth]{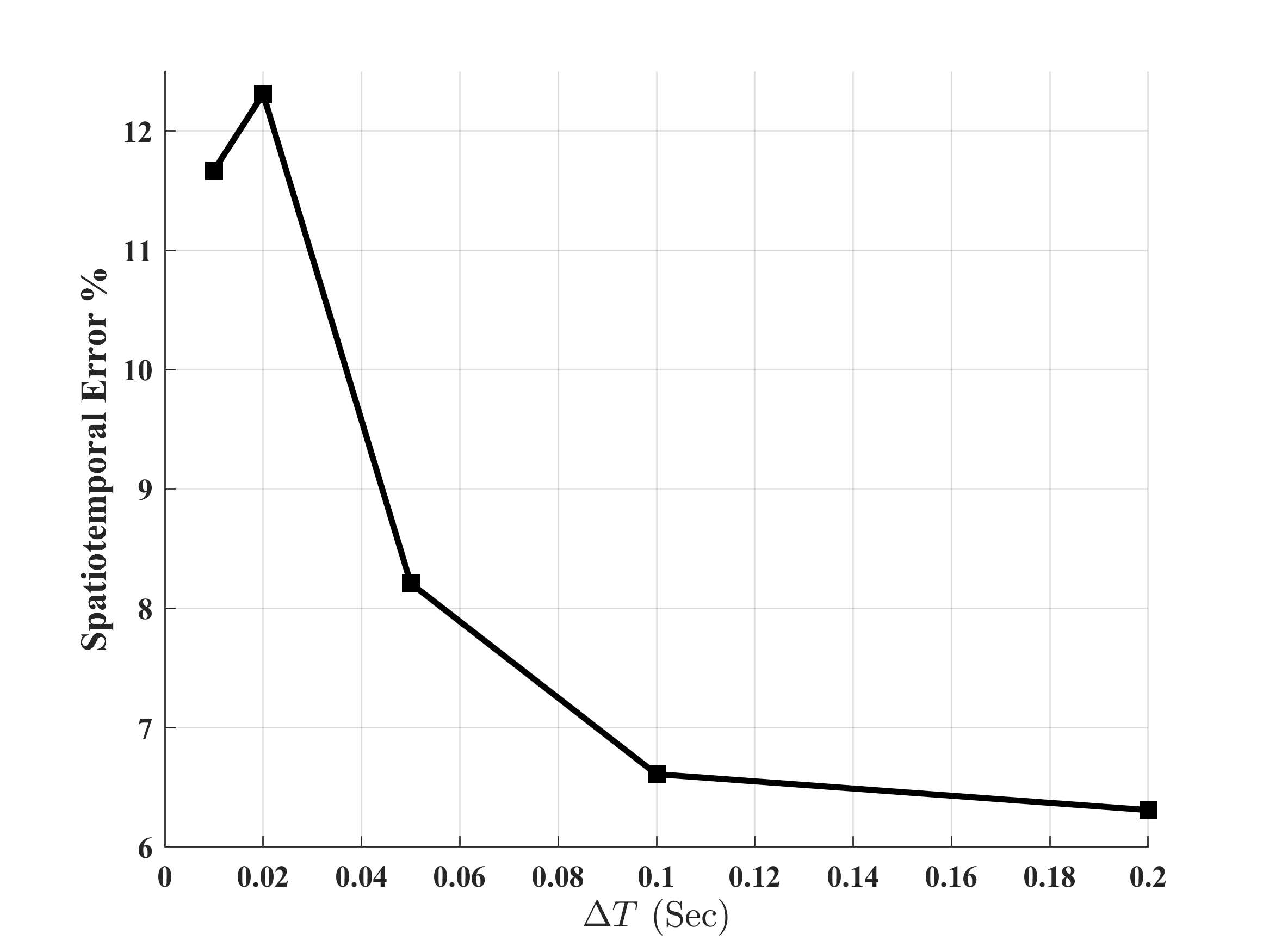}
        \caption{}
    \end{subfigure}
    \caption{Impact of $\Delta$t on the spatiotemporal error; (a) Gaussian (b) Multiquadric.}
    \label{Fig:DeltaT}
\end{figure*}

The approach also has strong sensitivity to time step changes as depicted in Figure \ref{Fig:DeltaT}. The error has a decreasing trend with increasing time steps for both Gaussian and Multiquadric kernels when the observation span remains constant.  It appears that the optimal time step for Multiquadric is 0.2 seconds, with a corresponding error of 7.31 percent. However, the minimum error for the Gaussian basis function is 7.59 percent when the time step is 0.1 second.\\

\begin{figure*}[h]
    \centering
    
    \begin{subfigure}{0.48\textwidth}
    \includegraphics[width=\linewidth]{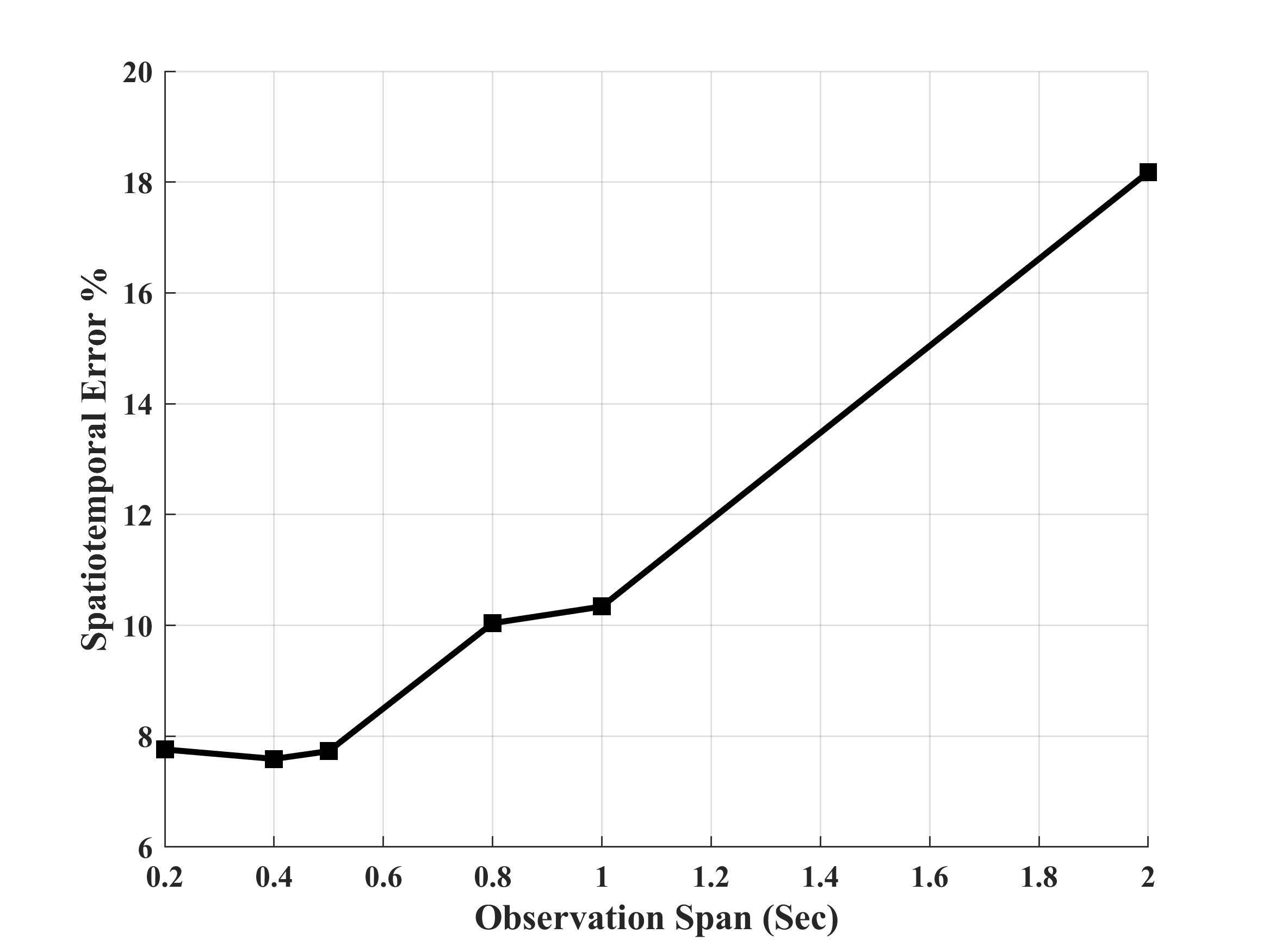} 
        \caption{}
    \end{subfigure}
    \hfill
    \begin{subfigure}{0.48\textwidth}
        \includegraphics[width=\linewidth]{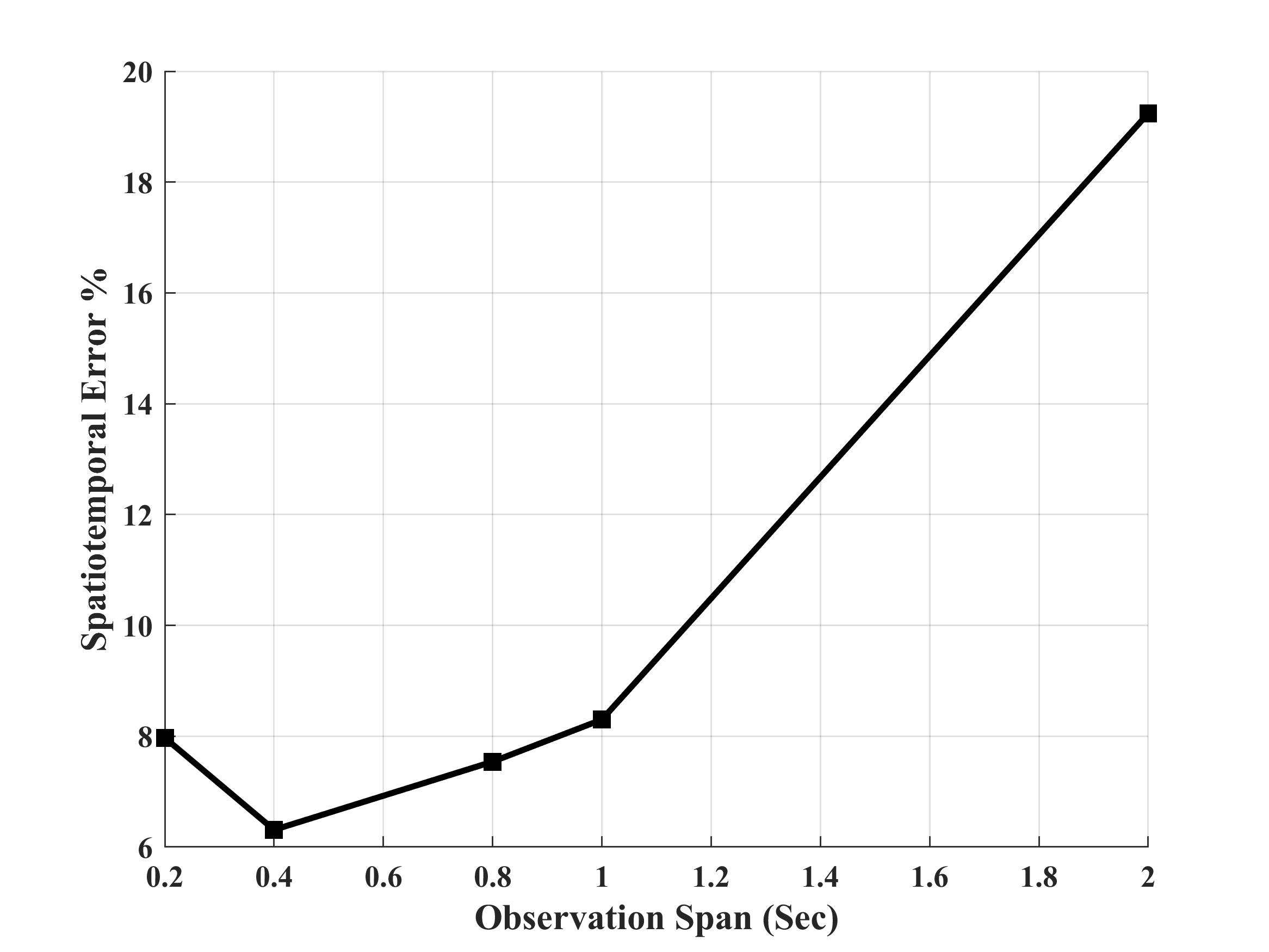}
        \caption{}
    \end{subfigure}
    \caption{Impact of observation span on the spatiotemporal error; (a) Gaussian (b) Multiquadric.}
    \label{Fig:Observation}
\end{figure*}

The duration over which temperature measurements can be obtained from thermocouples is a crucial factor in DA methods. Therefore, the length of the observation span can impact the accuracy of the proposed approach in estimating HF as parameters and temperature as states. As demonstrated in the line graph \ref{Fig:Observation}, a longer span leads to an increasing error for both Kernels. At $0.40$ seconds, the minimum error occurs for both kernels.\\

\begin{figure}[h]
    \centering
    
    \begin{subfigure}{0.7\textwidth}
        \includegraphics[width=\linewidth]
        {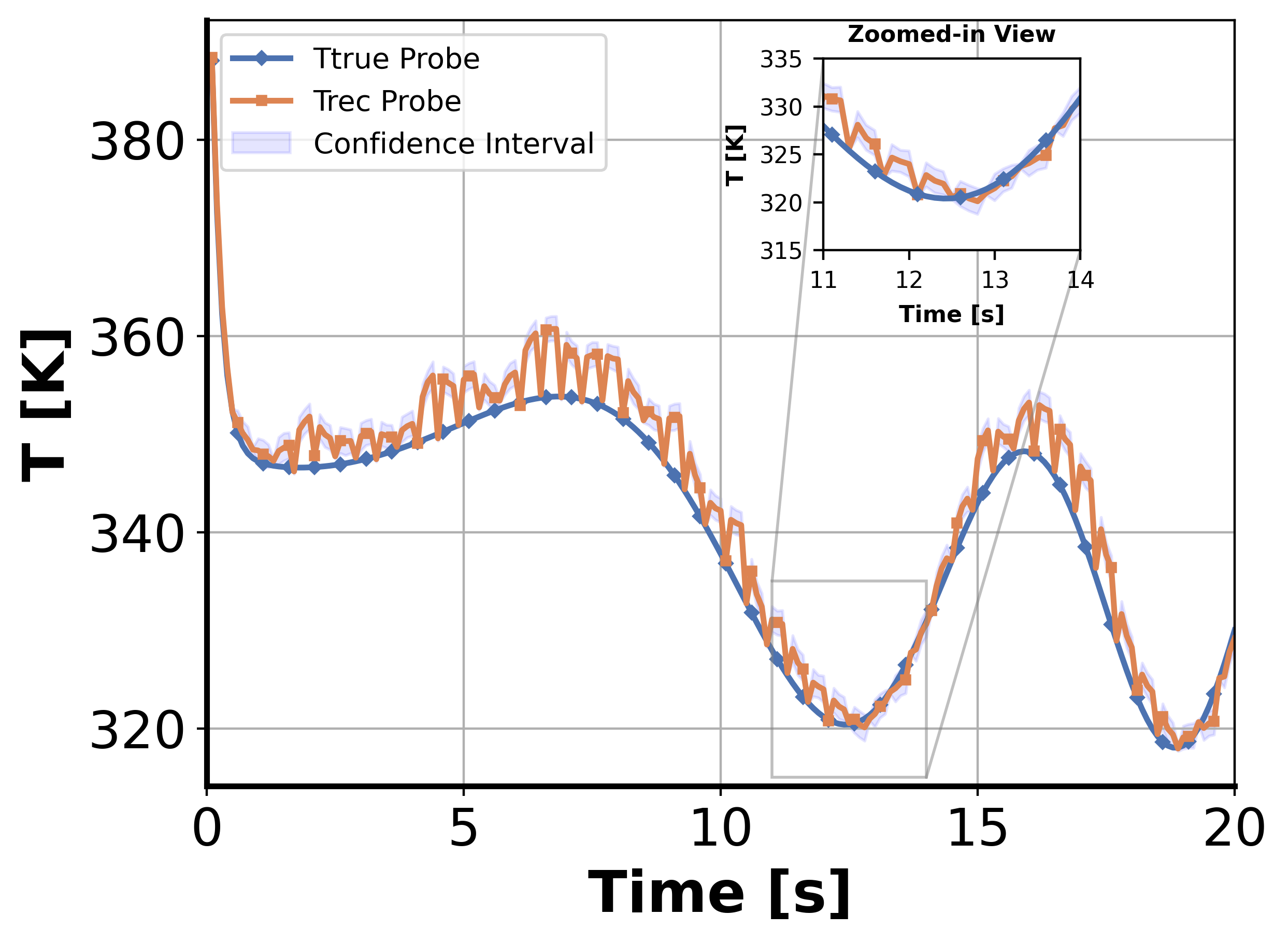}
        \caption{}
        \label{fig:TrueAndReconstructedMeanTemperatureAtaProbe_0.91_0.02_0.55_OverTimeGaussian}
    \end{subfigure}
    
    \begin{subfigure}{0.7\textwidth}
        \includegraphics[width=\linewidth]
        {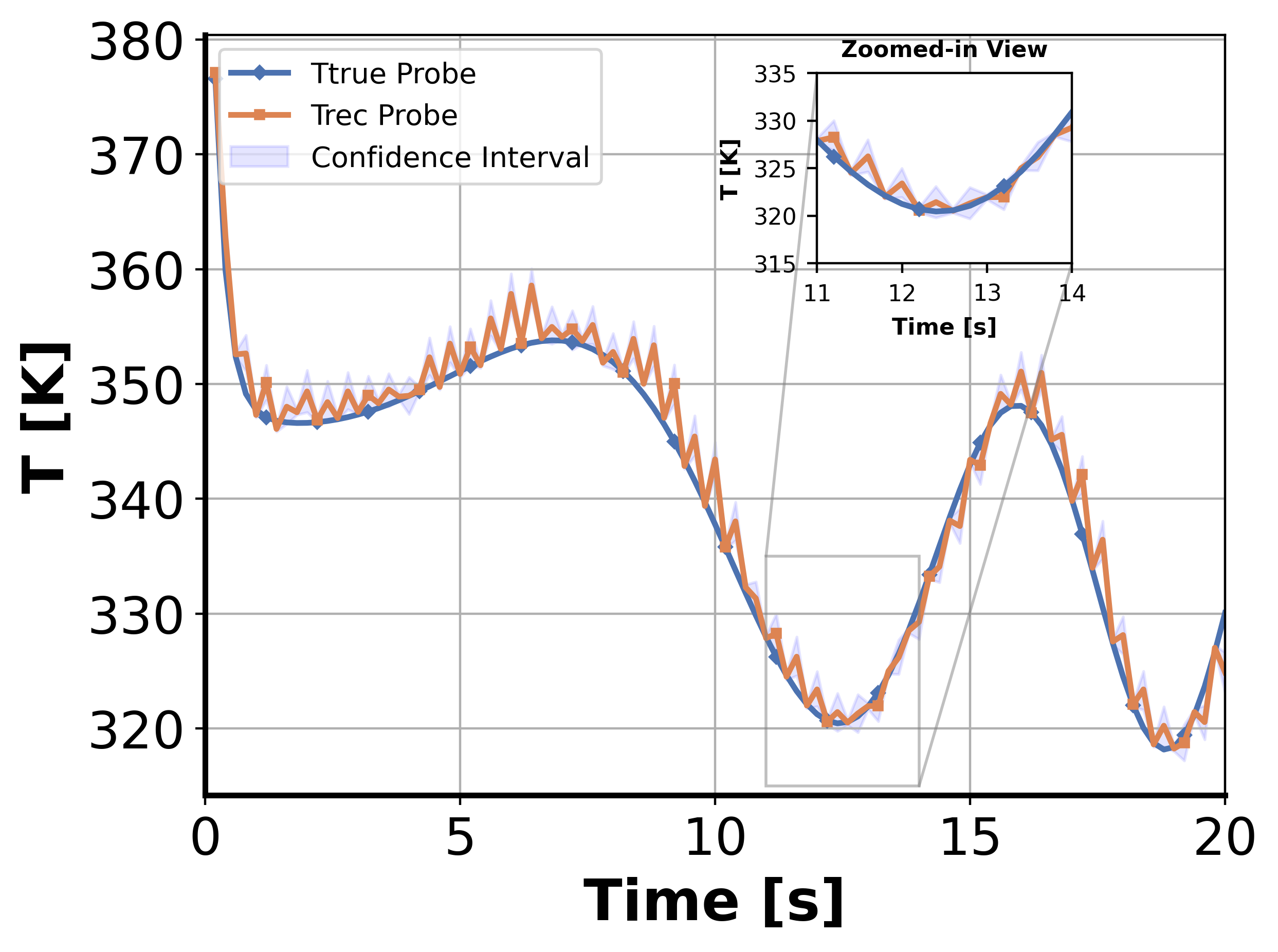}
        \caption{}
        \label{fig:TrueAndReconstructedMeanTemperatureAtaProbe_0.91_0.02_0.55_OverTimeMultiquadric}
    \end{subfigure}
    \caption{Comparison of true and predicted temperatures at the probe $(0.91, 0.02, 0.55)$ with confidence intervals (a) Gaussian and (b) Multiquadric kernels}
    \label{fig:probe_Temperature_comparison}
\end{figure}

To evaluate the accuracy and effectiveness of the proposed approach in estimating the state of Equations \eqref{eq:3DUnsteadyHeatconduction} and \eqref{eq:BoundaryCondition}, the graph presented in Figure \ref{fig:probe_Temperature_comparison} depicts the mean temperature predicted by the EnSISF-wDF (coral line) with the corresponding true temperature measured at a specific probe located at the point(0.91, 0.02, 0.55) (blue line) along with plotting the Ensemble taking the envelope of probability distribution between 5\% and 95\% for Gaussian and Mutiquadric kernels based on the optimal parameters shown in Table \ref{tab:Optimalparameters}. As previously mentioned, the precise location of each thermocouple is at the plane y=0.02 m, which is a few millimeters inward from the $\Gamma_{S_\text{in}}$. This temperature probe is the point that is also utilized as one of the locations of thermocouples. Therefore, the true temperature value stays within the confidence since the thermocouples are the only places in which the information is available. \\

\begin{table}[h]
\caption{Optimal parameters for EnSISSF incorporating Gaussian and Multiquadric kernels.}
\label{tab:Optimalparameters}
\centering
\begin{tabular}{|c c c c c c c c|}
\hline

\rowcolor{lightgray} 
\cellcolor[HTML]{FFFFFF} & \bm{$S_{n}$} & \bm{$\eta$} & \bm{$\kappa$} & \textbf{Prior weight shifting} & \bm{$\Delta t$} & \textbf{Obervation span} & \textbf{Error} \\ 
\hhline{~|}

\cellcolor{lightgray}\textbf{Gaussian}& \cellcolor[HTML]{FFFFFF}{\color[HTML]{343434} 375} & 0.5& 0.2 & 0 & 0.1  & 0.4 & 7.59\\ 
\hhline{~|}

\cellcolor{lightgray}\textbf{Multiquadric}
& 300 & 3 & 0.2 & 0.3  & 0.2 & 0.4 & 6.31
\\ \hline

\end{tabular}
\end{table}

\begin{figure}[h]
    \centering
    
    \begin{subfigure}{0.7\textwidth}
        \includegraphics[width=\linewidth]{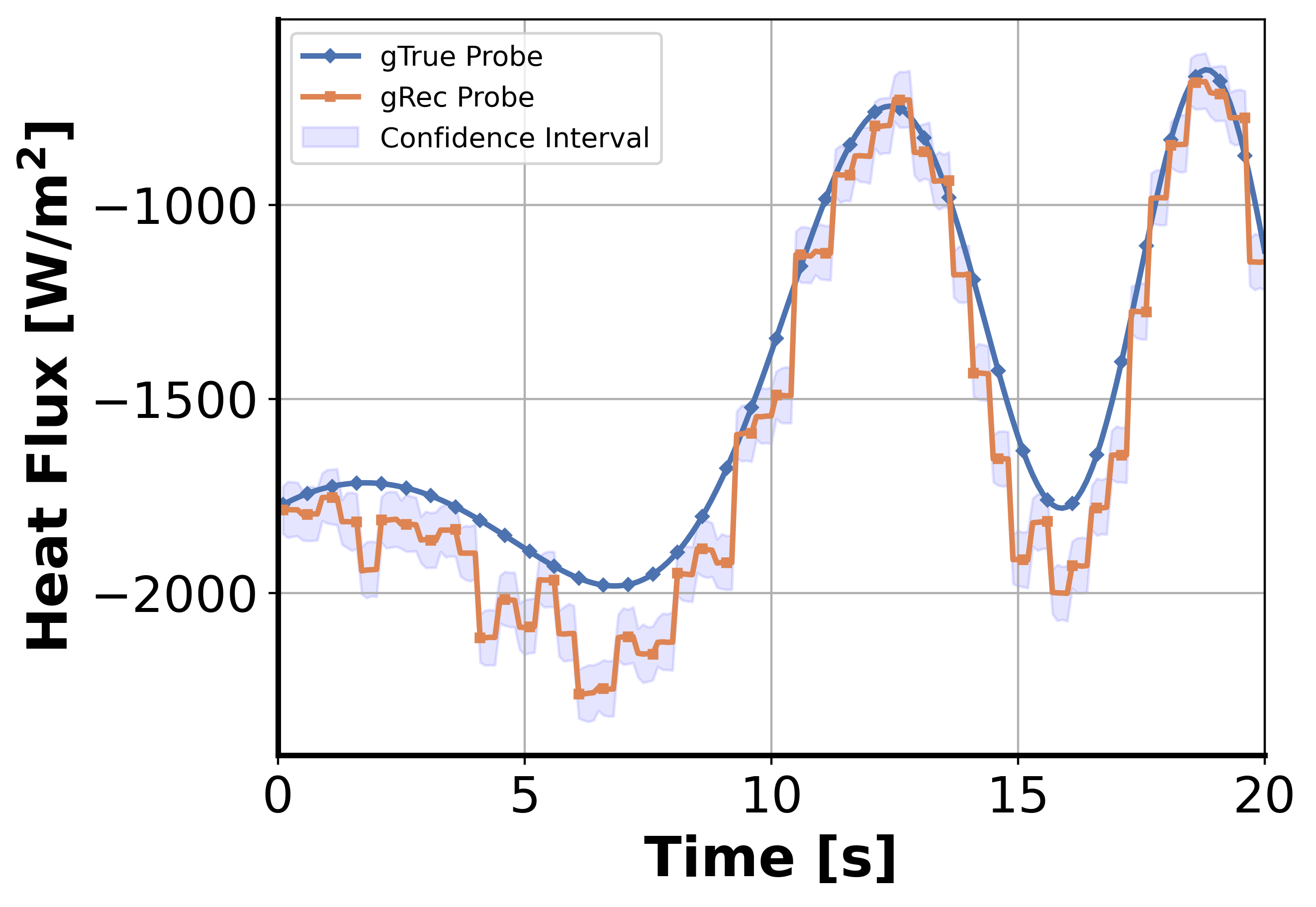}
        \caption{}
        \label{fig:TrueAndReconstructedMeanHeatFluxAtProbe_0.91_0.0_0.55_over_timeGaussian}
    \end{subfigure}
    
    \begin{subfigure}{0.7\textwidth}
        \includegraphics[width=\linewidth]{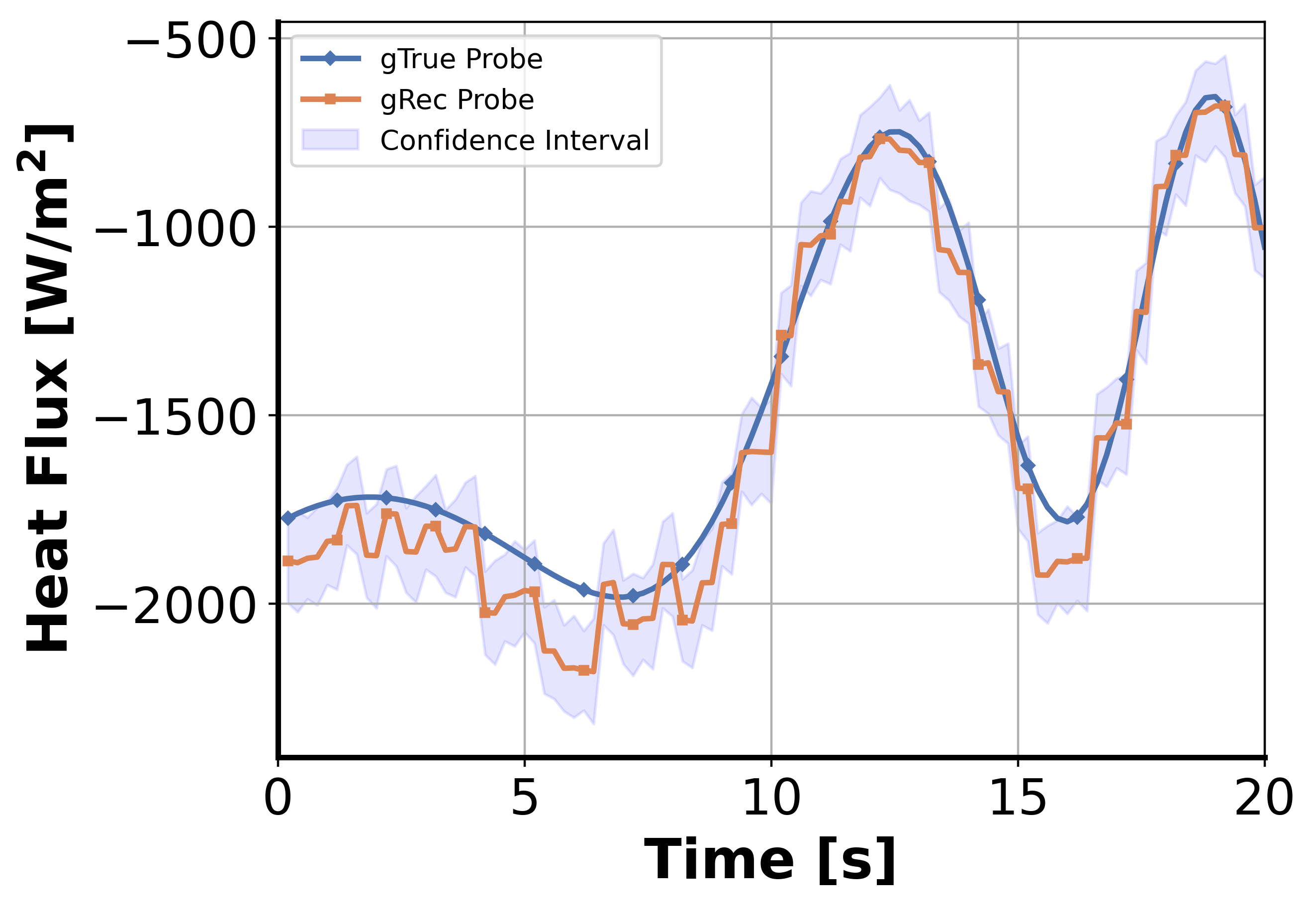} 
        \caption{}
        \label{fig:TrueAndReconstructedMeanHeatFluxAtProbe_0.91_0.0_0.55_over_timeMultiquadric}
    \end{subfigure}
    \caption{Comparison of actual and estimated HF at the probe (0.91, 0.00, 0.55) with
confidence intervals (a) Gaussian and (b) Multiquadric Kernels.}
    \label{fig:probe_HF_comparison}
\end{figure}

Figure \ref{fig:probe_HF_comparison} depicts a comparison between actual and estimated HF for the probe $(0.91, 0.0, 0.55)$ at $\Gamma_{S_\text{in}}$ with confidence intervals for Gaussian, and Multiquadric kernels. The normal vector of hot side BC goes from the $\Gamma_{S_\text{in}}$ as illustrated in Figure \ref{fig:Mold} to the molten steel. Thus, negative values of the HF over time mean that the flux goes opposite to the normal vector of $\Gamma_{S_\text{in}}$ going from outside to the inside.\\

\begin{figure}[h]
    \centering
    
    \begin{subfigure}{0.7\textwidth}
        \includegraphics[width=\linewidth]{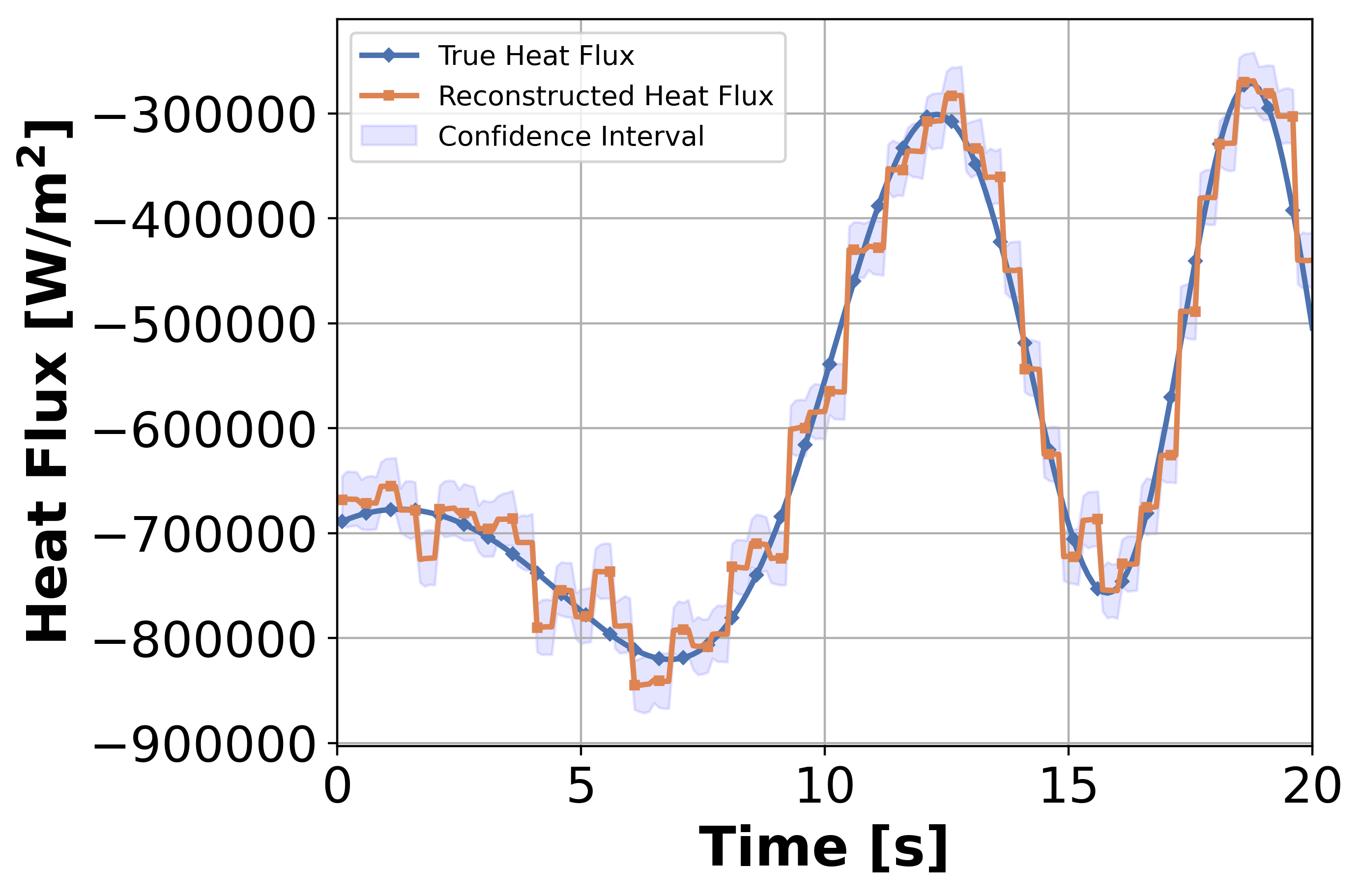}
        \caption{}
        \label{fig:TrueAndReconstructedMeanHeatFluxAtTheHotSideWithConfidenceIntervalGaussian}
    \end{subfigure}
    
    \begin{subfigure}{0.7\textwidth}
        \includegraphics[width=\linewidth]{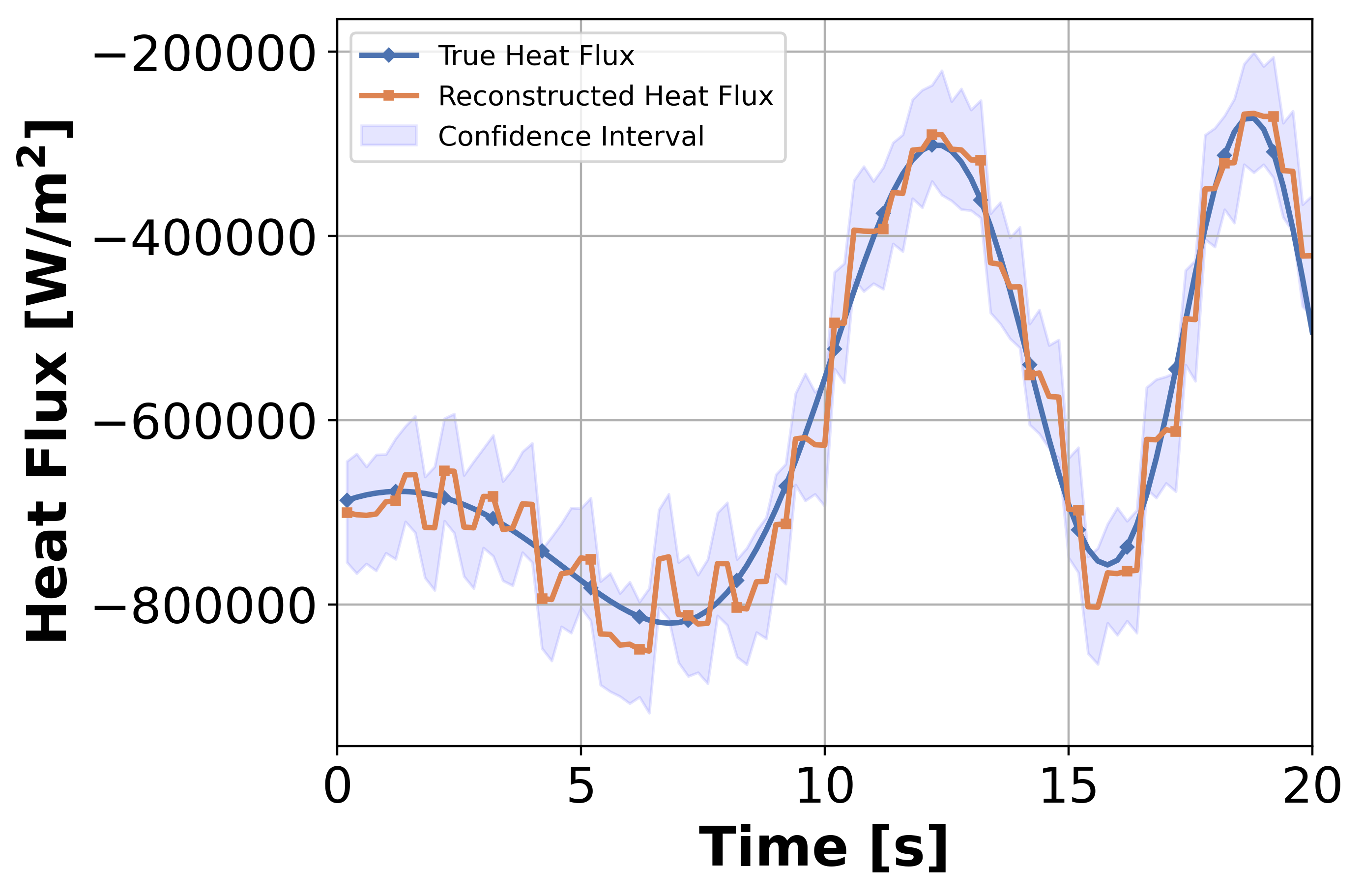} 
        \caption{}
        \label{fig:TrueAndReconstructedMeanHeatFluxAtTheHotSideWithConfidenceIntervalMultiquadric}
    \end{subfigure}
    \caption{Comparison of actual and estimated HF at entire $\Gamma_{S_\text{in}}$ with
confidence intervals (a) Gaussian and (b) Multiquadric Kernels.}
    \label{fig:probe_Total_HF_comparison}
\end{figure}

For each time interval in which an update phase is performed, additional measurements are collected. Thus, the confidence level is squeezed, leading to a reduction in the confidence level. Furthermore, Figure \ref{fig:probe_Total_HF_comparison} demonstrates the summation of HF values for all faces of $\Gamma_{S_\text{in}}$ over time which provides the total HF across the entire $\Gamma_{S_\text{in}}$ for two different kernels.\\ 

\begin{figure}[h]
    \centering
\includegraphics[scale=0.5]{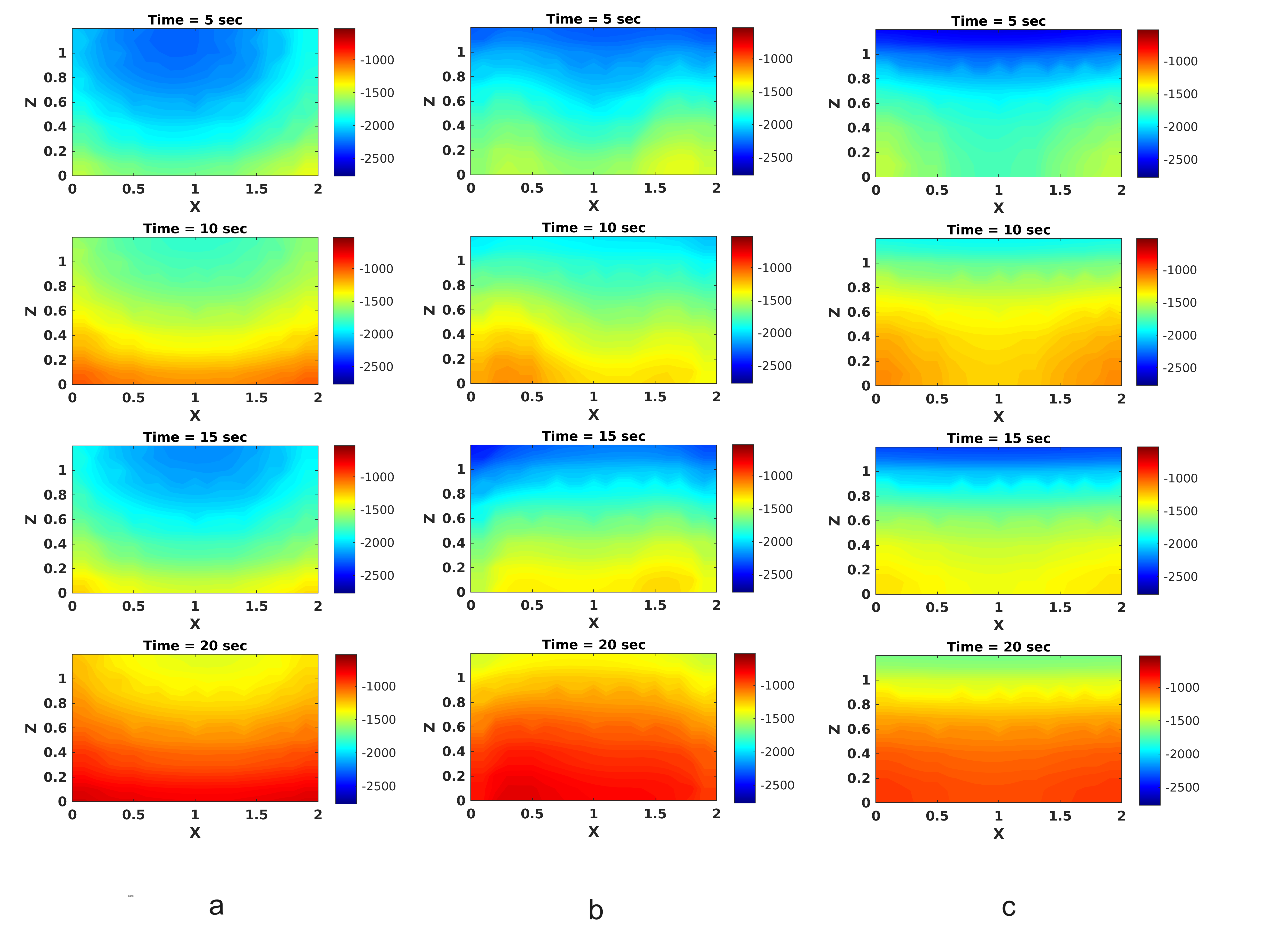}
    \caption{Contours of predicted HF at $\Gamma_{S_\text{in}}$ (a) Gaussian RBF (b) Multiquadric RBF (c) True HF.}
    \label{fig:SnapshotCountors}
\end{figure}

The HF contours presented in Figure \ref{fig:SnapshotCountors} offer a visual understanding of the predictive capabilities of the EnSISF with Gaussian and Multiquadric RBFs alongside the true HF at $\Gamma_{S_\text{in}}$ at various time steps (5, 10, 15, and 20 seconds). Upon visual inspection of the contours in Figures \ref{fig:SnapshotCountors}a and \ref{fig:SnapshotCountors}b, it is evident that the Multiquadric RBF produces contours that closely resemble the true HF contours (Figure \ref{fig:SnapshotCountors}c). This suggests that the Multiquadric RBF captures the spatial distribution of the HF at the boundary more accurately compared to the Gaussian RBF. To quantify the disparities between the predicted and true HFs, we computed the spatiotemporal relative error for each RBF as shown in Table \ref{tab:Optimalparameters}. The error for the Multiquadric RBF was found to be 6.31, while the Gaussian RBF yielded a slightly higher error of 7.59. This numerical assessment further supports the visual observation, indicating that the Multiquadric RBF provides a more accurate prediction of the true HF shown in Equations \eqref{eq:gTrue}, \eqref{eq:ConstantforgTrue} at $\Gamma_{S_\text{in}}$ boundary.\\

A further interesting investigation arises from the consideration of different true HF values. Since our study focused on a specific true HF scenario, it prompts us to contemplate whether the observed performance differences between RBFs hold under varying true HF conditions. It may be beneficial to explore alternative true HF scenarios and assess whether another type of RBF could outperform both Gaussian and Multiquadric RBFs under different system dynamics.\\

%
\newpage \clearpage
\section{Conclusion}\label{sec:conc}

The aim of the study was to employ a Bayesian approach called EnSISF-wDF incorporated with RBFs not only to predict the temperature distribution but also to estimate the HF between mold and molten steel by having some noisy temperature measurements inside the mold and considering process noise. By obtaining the probability distribution of HF with a certain mean and a minimal variance and comparing it with the true heat flux, it was proved that the EnSISF-wDF incorporating RBFs was reliable. However, solving IHTP in a deterministic framework only offered a singular mean value solution for the HF. The deterministic outcome failed to provide insights into the variability of the HF. Additionally, the procedure's dependence on certain hyperparameters was investigated that has not been addressed in the literature: the number of seeds, shape parameter, time step, observation window, and shifting mean of the prior weights of the RBFs and their covariance. Moreover, we proved that using EnSISF-wDF incorporating the Multiquadric kernel is more accurate compared to the EnSISF-wDF with Gaussian kernels. It was also shown that the former one requires less computational cost due to the smaller number of seeds that it needs compared to the latter one in order to reconstruct the HF for this case study. \\

Further investigation into the utilization of whether an intrusive reduced-order modeling strategy or surrogate models for real industrial problems would be advantageous in the future as well, as it would help alleviate the computational burden linked to the dynamic model of the forward problem of the algorithm. Moreover, investigating different RBFs to variations in true HF scenarios guides the selection of appropriate RBFs for specific operating conditions. Future research in this direction will ensure the versatility of predictive models, making them adaptable to a range of scenarios encountered in the CC process. Investigating different numbers of RBFs to find the optimal one depending on the complexity of the HF is another hyperparameter that should be studied.\\

\newpage \clearpage
\section*{Acknowledgements}
This work was supported by PRIN “Numerical Analysis for Full and Reduced Order
Methods for the efficient and
accurate solution of complex systems governed by Partial Differential Equations” (NA-FROM-PDEs) project.\\

\bibliographystyle{unsrt}
\bibliography{biblio.bib}
\end{document}